\documentclass[preprint,12pt]{elsarticle}
\biboptions{numbers,sort&compress}

\usepackage{amsthm}
\usepackage{amsmath,amssymb,amsfonts}
\usepackage{algorithmic}
\usepackage{graphicx}
\usepackage{textcomp}

\newtheorem{theorem}{Theorem}
\newtheorem{remark}{Remark}
\newtheorem{lemma}{Lemma}
\newtheorem{proposition}{Proposition}
\newtheorem{definition}{Definition}
\usepackage{appendix}
\usepackage{setspace}
\usepackage{multirow}
\usepackage{xcolor}
\usepackage{bbding}
\usepackage{booktabs}
\usepackage{algorithm}
\usepackage{caption}
\usepackage{subcaption}
\usepackage{threeparttable}
\usepackage{hyperref}
\newtheorem{problem}{Problem}

\journal{European Journal of Control}

\begin{document}

\begin{frontmatter}



\title{A Time-Triggered Dimension Reduction Algorithm \\ for the Task Assignment Problem}

\author{Han Wang}\ead{han.wang@eng.ox.ac.uk}
\author{Kostas Margellos}\ead{kostas.margellos@eng.ox.ac.uk}
\author{Antonis Papachristodoulou}\ead{antonis@eng.ox.ac.uk}

\address{Department of Engineering Science, University of Oxford, Parks Road, Oxford, OX1 3PJ, United Kingdom}

\begin{abstract}
The task assignment problem is fundamental in combinatorial optimisation, aiming at allocating one or more tasks to a number of agents while minimizing the total cost or maximizing the overall assignment benefit. This problem is known to be computationally hard since it is usually formulated as a mixed-integer programming problem. In this paper, we consider a novel time-triggered dimension reduction algorithm (TTDRA). We propose convexification approaches to convexify both the constraints and the cost function for the general non-convex assignment problem. The computational speed is accelerated via our time-triggered dimension reduction scheme, where the triggering condition is designed based on the optimality tolerance and the convexity of the cost function. Optimality and computational efficiency are verified via numerical simulations on benchmark examples.

\end{abstract}




\end{frontmatter}


\section{Introduction}
The Task Assignment Problem (TAP) is of great importance in combinatorial optimisation \cite{papadimitriou1998combinatorial}. The formulation of this problem is that, given a set of agents and a set of tasks, each agent can select a task from its admissible task set, while pairing between tasks and agents is one-to-one. The goal is to minimize the total cost or maximize the global utility. 

Typical types of TAPs are the linear assignment problem (LAP) and the quadratic assignment problem (QAP). The LAP has a linear cost function and binary decision variables, i.e., $1$ means a selection and $0$ means a rejection. There are many real-world applications that can be formulated as a LAP, ranging from robot formation control \cite{michael2008distributed}, \cite{macdonald2011multi}, \cite{morgan2016swarm} to facilities allocation \cite{hillier1966quadratic}. The LAP has been well investigated with many different methods in both centralized and decentralized ways. Centralized algorithms include the Hungarian algorithm \cite{kuhn1955hungarian}, and its extensions involves iterative algorithms \cite{bertsekas1981new}, \cite{hung1983polynomial}. These methods have high computational efficiency for small-scale problems, but they cannot be easily implemented in parallel, limiting their usage to large-scale instances. To address this limitation, another class of algorithms called \emph{auction algorithms} \cite{zavlanos2008distributed}, \cite{choi2009consensus}, \cite{bertsekas2009auction}, \cite{leyton2000towards} provide an improved solution, imitating bidding in an auction. For a comprehensive literature review of this area, readers are referred to the survey paper \cite{khamis2015multi}.

Another type of TAP is the QAP, which was firstly proposed in \cite{koopmans1957assignment} aiming at solving resource allocation problems. The QAP has a similar structure to the LAP, except for the quadratic cost function. QAPs are widely used in many applications, e.g. the travelling salesman problem \cite{applegate2011traveling}, graph matching \cite{williams1997multiple}, \cite{zhou2015factorized}, \cite{cour2007balanced} etc. Unlike the LAP, which allows finding a global optimum in polynomial time with multiple heuristic algorithms, the QAP is provably NP-hard in general~\cite{papadimitriou1998combinatorial}, which makes it hard to design efficient heuristic
algorithms. Relaxation algorithms have been proposed to overcome this limitation. The original mixed-integer programming problem is relaxed into a continuous optimisation problem \cite{zhao1998semidefinite}. This idea is derived from the fact that permutation matrices are at the vertices of Birkhoff polytopes, i.e., the class of doubly stochastic matrices. Thus, several efforts have been made to find good solutions over doubly stochastic matrices, followed by a projection of the solution to the set of permutation matrices. A penalty term can be added to the original quadratic cost function to convexify the nonconvex quadratic cost function \cite{fogel2013convex}. Tighter convex underapproximations were proposed to improve the optimality of the solution \cite{dym2017ds++} \cite{bernard2018ds}.  

The projection on the space of permutation matrices is realized with a convex-to-concave method \cite{zaslavskiy2008path}, which incrementally tunes the penalty quantity to change the convexity of the quadratic cost from convex to concave. The solution is a permutation matrix since the minimizer lies on the vertices of the convex constraint hull while minimizing a concave cost function \cite{boyd2004convex}. This algorithm finds a solution path of a family of convex to concave minimization problems, obtained by linearly interpolating between the convex and concave relaxations. The interpolation procedure requires solving additional $n$ problems, which makes it computationally expensive. Besides, the sub-problems in the procedure are nonconvex.

Motivated by the above limitations, in this paper, we propose a fast relaxation-based iterative algorithm for the TAP, especially for the nonconvex QAP. Our contributions are twofold:
\begin{itemize}
    \item A convex relaxation-based framework is proposed. The non-convex permutation set is relaxed to a polyhedral doubly stochastic set. The nonconvex quadratic cost function is convexified to be $\sigma$-strongly convex.
    \item Computational speed acceleration is realized through a time-triggered dimension reduction approach: we reduce the dimension of the decision matrix incrementally by removing columns and rows. Compared to the convex-to-concave method, our method only solves the optimisation problem once, and the dimension of the problem is incrementally decreased. In addition, the convexity of the convexified cost function is preserved across iterations.
\end{itemize}

The remainder of the paper is organized as follows. In Section \ref{sec:2} we define the notations used in this paper and introduce the formulation of TAP, including LAP and QAP. The convex relaxation is presented in Section \ref{sec:4}. Section \ref{sec:5} details our time-triggered dimension reduction algorithm, followed by the simulations in Section \ref{sec:6}. Section \ref{sec:conclusion} concludes the paper.

\section{Problem Formulation}
\label{sec:2}

\subsection{Preliminaries}
We use $I_N$ to denote the $N \times N$ identity matrix, and ${\bf 1}_N$ to denote the $N$-dimensional vector whose entries are all equal to one. When it is obvious from the context, we omit the subscripts and use $I$ and $\bf{1}$ instead, respectively. $\mathbb{R}$ represents the set of real numbers and $|| \cdot ||_F$ denotes the Frobenius norm. We let $\mathrm{vec}(X)$ denote the vectorizartion operation for a matrix $X$, and $\mathrm{vec}^{-1}(x)$ the inverse operation that takes as input a column vector and returns a matrix for a vector $x$. $A \otimes B$ represents the Kronecker product of matrix $A$ and $B$. $X\ge0$ denotes point-wise non-negativity of the elements of matrix $X$. $X \succeq 0$ means that matrix $X$ is positive semi-definite. Let $\mathbb{DS}_N$ denote the set of $N \times N$ doubly stochastic matrices, i.e. $\mathbb{DS}_N=\{X: X\ge0,\textbf{1}^TX=\textbf{1}^T,X\textbf{1}=\textbf{1}\}$, $\Pi_N$ denote the set of $N \times N$ permutation matrices, i.e. $\Pi_N=\{X\in\{0,1\}^{N\times N}:X^TX=I_N\}$. $\text{Tr}(X)$ denotes the trace of matrix $X$. $\mathrm{th}_{\ge 0}(X)$ sets the negative elements of $X$ to zero.

\begin{definition}\label{def:convex}
A differentiable function $f(\cdot)$ is called \textit{$\sigma$-strongly convex} with $\sigma>0$ on a domain $\mathcal{D}$ if the following inequality holds for all $x,y \in \mathcal{D}$:
\[f(y) \ge f(x) + \nabla f{(x)^T}(y - x) + \frac{\sigma }{2}||y - x|{|^2_2}.\]
\end{definition}

We now discuss the LAP and QAP. Although the LAP can be regarded as a special case of QAP, we still want to briefly introduce it for its importance in applications.

\subsection{Linear Assignment Problem}
 The LAP describes a scenario in which every agent $i$ is capable of choosing one particular task $j$ from the tasks pool. After that, a specific predefined work cost or award is added. Each agent can only select a single task, and every task can only be allocated to one agent. Our goal is to find an optimal assignment strategy for each task/agent pair, to realize a maximum reward or minimum cost. 

In considering such a problem of pairing agents with tasks, a linear programming formulation can be derived
\begin{equation}\label{eq1}
\begin{array}{l}
\mathop {\min }\limits_X F(X)=\sum\limits_{i = 1}^N {\left( {\sum\limits_{j = 1}^N {{X_{i,j}}}\beta_{i,j}} \right)} \\
\mathrm{subject~to}~X \in \Pi{_N},\\
\end{array}
\end{equation}
where $\beta_{i,j}\in \mathbb{R}$ is the predefined cost for agent $i$ to select task $j$. Decision variable $X$ should be within the $N$-dimensional permutation set. The LAP is a tractable P-problem, which has an efficient solution without additional constraints.



\subsection{Quadratic Assignment Problem}
The Quadratic assignment problem originates from facility-location allocation problems. Suppose there are $n$ facilities and $n$ locations, and assume that the distances between locations and flows of facilities are known. The problem is to assign the facilities to different locations resulting in a minimum sum of the products of the distances and flows. Intuitively, pairwise facilities which have higher flows are encouraged to be placed at nearer locations. Unlike the linear assignment problem, the cost function of this problem is expressed in terms of a quadratic function. 

The classic Koopman-Beckmann formulation of the QAP is as follows
\begin{equation}\label{kb}
    \begin{array}{l}
\mathop{\min}\limits_X \text{Tr}(AXB{X^T})\\
\mathrm{subject~to}~X \in {\Pi _N},
\end{array}
\end{equation}
where $A \in \mathbb{R}^{N\times N}$ is the flows matrix and $B \in \mathbb{R}^{N\times N}$ is the distance matrix. Making use of the cyclic properties of the $\text{Tr}(\cdot)$, we can rewrite the objective function in \eqref{kb} as follows
\begin{equation}\label{tr}
        \begin{split}
\text{Tr}(XB{X^T}A) &= \text{Tr}({X^T}AXB)\\
 &= \mathrm{vec}{(X)^T}\mathrm{vec}(AXB)\\
 &= \mathrm{vec}{(X)^T}(B \otimes A)\mathrm{vec}(X)\\
 &= {x^T}Wx,
\end{split}
\end{equation}
where $x=\mathrm{vec}(X)$ is the vectorization of the permutation matrix $X$, $W\in \mathbb{R}^{N^2\times N^2}=(B \otimes A)$. In this paper we consider a more general form of \eqref{tr}, where a linear term is added:

\begin{equation}\tag{QAP}\label{OQAP}
\begin{array}{l}
\mathop {\min }\limits_x f(x)=xWx^T+c^Tx\\
\mathrm{subject~to}~\mathrm{vec}^{-1}(x) \in \Pi_N,\\
\end{array}
\end{equation}
where $W \in \mathbb{R}^{N^2 \times N^2}$, $c \in \mathbb{R}^{N^2}$. It should be noted that the LAP is a special case of the QAP with zero quadratic matrix $W$. In the sequel we just consider the QAP as it covers the LAP as a special case. Here we denote this problem as the original QAP.

Existing results point out that it is hard to solve \eqref{OQAP} precisely because it is NP-hard \cite{papadimitriou1998combinatorial}. Possible additional constraints can also increase the complexity, even when the constraints are affine. This motivates us to develop an efficient algorithm to find a sub-optimal solution for problem \eqref{OQAP}, as well as address the challenge of additional constraints.




\section{Convex Relaxation}
\label{sec:4}

To efficiently solve OQAP, an intuitive idea is to use constrained quadratic programming tools. However, \eqref{OQAP} is generally nonconvex for two reasons: i) the constraint is nonconvex; ii) the cost function is nonconvex when $W$ is not necessarily positive semi-definite. In this section we present results with respect to the convexification, including the convexification of the constraints $\mathrm{vec}^{-1}(x)\in\Pi_{N}$ and the cost function $f(x)$. We begin this section with the convexification of the constraint set, i.e. the set of permutation matrices.

\subsection{Convexifying the Constraints}
The original problem \eqref{OQAP} is a mixed-integer programming problem over permutation matrices. $X_{ij} = 1$ implies that a specific task is chosen by a corresponding agent, and \emph{vice-versa}. Intuitively, the resulting minimizer of such problems is a deterministic distribution, which is a special case of a stochastic distribution. This motivates us to consider solving a relaxed version of the original assignment problem over the set of \emph{doubly stochastic matrices}, which renders the optimisation problem continuous.  

The constrained relaxed QAP \eqref{CRQAP} is described as:
\begin{equation}\tag{CRQAP}\label{CRQAP}
    \begin{array}{l}
\mathop {\min }\limits_{ x} {{f}}( x)\\
\mathrm{subject~to}~\mathrm{vec}^{-1}( x) \in \mathbb{DS}_N.
\end{array}
\end{equation}

\subsection{Convexifying the Cost}
Convexifying the nonconvex quadratic cost function has been well studied in quadratically constrained quadratic programming (QCQP) over the past decades \cite{anstreicher2012convex}. Existing QAP papers \cite{fogel2013convex}, \cite{dym2017ds++}, \cite{bernard2018ds}, \cite{ferreira2018semidefinite} proposed a series of convexification tools which can be summarized into two categories, $\alpha$BB convexification $\widehat {QAP}$, and semidefinite programming (SDP) convexification $\widetilde {QAP}$. The first one is defined as follows:
\begin{equation}\tag{$\widehat {QAP}$}\label{abb}
\begin{split}
    &{x^T}(W + \mathrm{diag}(\alpha ))x + {(c - \alpha )^T}x,
\end{split}
\end{equation}
where $\alpha \in \mathbb{R}^{N}_+$ is chosen so that $W+\mathrm{diag}(\alpha)$ is positive semi-definite. \cite{fogel2013convex} selected $\alpha=-\mu_{\min}\bf{1}$, where $\mu_{\min}$ denotes the minimum eigenvalue of $W$; \eqref{abb} is therefore convex over $\mathbb{R}^N$. \cite{bernard2018ds} further proposed a tighter convexification by restricting the variable to an affine space $X=X_0+Fz$, where $F\in\mathbb{R}^{N^2\times (N-1)^2}$ denotes the null space of $\mathbb{DS}_N$.   

The SDP convexification $\widetilde {QAP}$ comes from an observation that the cost function in \eqref{tr} could be rewritten as:
\[\text{Tr}((B \otimes A)\mathrm{vec}(X)\mathrm{vec}{(X)^T}).\]
Substituting $Q=\mathrm{vec}(X)\mathrm{vec}{(X)^T}$, \eqref{OQAP} becomes:
\begin{equation*}
    \begin{split}
\mathop {\min }\limits_{Q,x} ~~~&\text{Tr}(WQ)\\
\mathrm{subject~to}~~~&Q = x{x^T}.
    \end{split}
\end{equation*}
Note that the cost function above is linear, the only nonconvex part is the constraint on $Q$. Such constraint can be relaxed to $Q-xx^T \succeq 0$, and by some means of the Schur complement this is equivalent to:
\[\left[ {\begin{array}{*{20}{c}}
Q&x\\
{{x^T}}&1
\end{array}} \right]\succeq 0.\]
Substituting the quadratic constraint $Q=xx^T$ with the above positive semi-definite constraint, the SDP convexification is formulated as:
\begin{equation}\tag{$\widetilde {QAP}$}
    \begin{split}
\mathop {\min }\limits_{Q,x}~~~& WQ\\
\mathrm{subject~to}~~~&\left[ {\begin{array}{*{20}{c}}
Q&x\\
{{x^T}}&1
\end{array}} \right]\succeq 0,\\
&x \ge 0.\\
    \end{split}
\end{equation}

\begin{theorem}\label{th1}
(\cite{anstreicher2012convex}) Let $\hat z$ and $\tilde z$ denote the optimal values of the problems $\widehat{QAP}$ and $\widetilde{QAP}$, respectively. Then $\hat z \le \tilde z$.
\end{theorem}

Theorem \ref{th1} shows that $\widetilde{QAP}$ gives a tighter underestimation than that of $\widehat{QAP}$ at the price of lifting the variable dimension, as well as solving a more complex SDP rather than a simpler QP. 

In this paper, we convexify the cost function based on $\widehat {QAP}$ to reach a trade-off between computational efficiency and tight underestimation. To obtain a better convergence rate, we convexify the cost function to be $\sigma$-strongly convex.

\begin{lemma}\label{lem:convex}
The cost function $f$ is $\sigma$-convex if and only if $W-\frac{\sigma}{2} I \succeq 0$ .
\end{lemma}

\begin{proof}
From the definition of $\sigma$-strong convexity we obtain, $\forall x,y$, $\mathrm{vec}^{-1}(x)\in \mathbb{DS}_{N}$, $\mathrm{vec}^{-1}(y)\in \mathbb{DS}_{N}$:
\begin{equation}\label{eq13}
f(\alpha x + (1 - \alpha )y) \le \alpha f(x) + (1 - \alpha )f(y) - \frac{{\alpha (1 - \alpha )}\sigma}{2}||x - y|{|^2_2},
\end{equation}
for any $\alpha\in [0,1]$, which implies that:
\begin{equation}\label{eq14}
\begin{split}
    &\alpha (1 - \alpha ){x^T}Wx + \alpha (1 - \alpha ){y^T}Wy - \alpha (1 - \alpha ){x^T}Wy\\
    &- \alpha (1 - \alpha ){y^T}Wx - \frac{\alpha(1-\alpha)\sigma}{2}||x-y||_2^2 \ge 0\\
    &\Leftrightarrow (x-y)^TW(x-y)-\frac{\sigma}{2}(x-y)^T(x-y) \ge 0\\
    &\Leftrightarrow (x-y)^T(W-\frac{\sigma}{2}I)(x-y)\ge 0.
\end{split}
\end{equation}
The reverse is similar. Thus, $W - \frac{\sigma }{2}I \succeq 0 \Leftrightarrow $ $f(x)$ is $\sigma$-strongly convex.
\end{proof}

Following the results of Lemma \ref{lem:convex}, the relaxed QAP (RQAP) formulation is:
\begin{equation}\tag{RQAP}\label{RQAP}
    \begin{array}{l}
\mathop{\min}\limits_x~\tilde f(x) = {x^T}(W + (\frac{\sigma }{2} - {\mu _{\min}}){I_N})x + {c^T}x+\frac{- \frac{\sigma }{2} + {\mu _{\min}}}{N}\\
\mathrm{subject~to}~~~\mathrm{vec}^{-1}(x) \in \mathbb{DS}_N.
\end{array}
\end{equation}

It can be seen that $\tilde f(x) \le f(x)$, since $\forall  x \in \mathbb{DS}_N,(\frac{\sigma}{2}-\mu_{\min})x^TI_Nx\le(\frac{\sigma}{2}-\mu_{\min})/N$. The equality is fulfilled only when $x\in \Pi_N$. In the sequel we use $\tilde W=W+(\frac{\sigma}{2}-\mu_{\min})I_N$ for brevity.

After presenting the convexification for both the constraint and the cost function, the comparison between the optimal values of the original problem \eqref{OQAP} and the relaxed problem \eqref{RQAP} is given by the following lemma. 
\begin{lemma}\label{lem:bound}
The optimal values of the relaxed problem \eqref{RQAP} and the original problem \eqref{OQAP} satisfy:
\begin{equation*}
    \mathop {\min }\limits_{\mathrm{vec}^{-1}(x) \in \mathbb{DS}_N} \tilde f(x) \le \mathop {\min }\limits_{\mathrm{vec}^{-1}(x) \in {\Pi _N}} f(x).
    \vspace{5pt}
\end{equation*}
\end{lemma}

\begin{proof}
According to the Birkhoff–von Neumann theorem, the set of $N \times N$ doubly stochastic matrices $\mathbb{DS}_N$ is the convex hull of the set of $N \times N$ permutation matrices $\Pi_N$, i.e. $\Pi_N \subset \mathbb{DS}_N$. Besides, we have $\tilde f(x)\le f(x)$ according to the previous discussion. Thus, the optimal value over the set of doubly stochastic matrices is no bigger than that over the set of permutation matrices.
\end{proof}

Lemma~\ref{lem:bound} reveals that the minimum of the relaxed problem \eqref{RQAP} is a lower bound of the minimum of the original problem \eqref{OQAP}. This immediately implies the following proposition about when a zero gap is realized. 

\begin{proposition}\label{pr1}
Let $x^*$ be the minimizer of \eqref{RQAP}. Then $x^*$ is the global minimizer of \eqref{OQAP} if and only if $\mathrm{vec}^{-1}(x^*) \in \Pi_{N}$.
\end{proposition}
Proposition \ref{pr1} shows that if the optimal solution of \eqref{RQAP} lies in the permutation set, we can conclude that we have found the global minimizer of the original QAP \eqref{OQAP}. 

\section{Dimension Reduction Algorithm}
\label{sec:5}

The relaxed problem \eqref{RQAP} is a smooth continuous optimisation problem, which can be solved with multiple numerical solvers, e.g. Ipopt, SQP. However, when the problem's dimension becomes exceptionally high, such a problem is hard to solve in real time. More specifically, consider the quadratic assignment problem: the dimension of variable $x$ is $N^2$, where $N$ is system size.  
\subsection{Time Triggered Dimension Reduction}
Here we notice that the permutation matrix has a sparsity structure, where only $N$ entries are non-zero out of $N^2$ elements. This motivates us to propose a dimension reduction scheme to accelerate the solution. For example, when one element increases faster than others, we can binarize it into one and then ground the elements which lie within the same column or row to zero. This also enables us to transform the decision variable from a doubly stochastic matrix to a permutation matrix. Our method is based on the steepest projection gradient method; we chose such a method compared to other first-order optimisation methods because the original problem has an explicit gradient expression at every iteration. Besides, with our alternating directional projection algorithm, the projection can be realized within a few steps. 

The time-triggered dimension reduction algorithm is described in Algorithm \ref{al1}.
\begin{algorithm}[h]
  \caption{Time-Triggered Dimension Reduction Algorithm (TTDRA)}
   \hspace*{\algorithmicindent} \textbf{Input:} parameters $\epsilon\in \mathbb{R}^+,\{\alpha_k\} \subseteq (0, + \infty )$\\
      \hspace*{\algorithmicindent} \textbf{Initialize}: $x(0)\in\mathbb{R}^{N^2},\tilde W(0)=\tilde W,l=0,\mathrm{count}=0,n_{\mathrm{iter}}$, $P=0$, $\eta$\\
 \hspace*{\algorithmicindent} \textbf{Output:} resulting permutation matrix $P$\\
 \vspace{-3ex}
  \begin{algorithmic}[1]\label{al1}
  \WHILE{$P\notin \Pi_N$}
  \STATE $x(k+1)\leftarrow Proj_{\mathbb{DS}}(x(k)-\alpha_k\nabla \tilde f(x(k)))$
  \IF {$\mathrm{count}\ge n_{\mathrm{iter}}$}
    \STATE $\mathrm{count} \leftarrow 1$
  \STATE $X(k+1) \leftarrow \mathrm{vec}^{-1}(x(k+1))$
  \STATE find index $\{c,r\}$ of the maximum element from $X(k+1)$
  \STATE delete column $c$ and row $r$ from $X(k+1)$
  \STATE $\tilde W(l+1) \leftarrow$delete columns and rows from $\tilde W(l)$
  \STATE find real original row and column index $\tilde c,\tilde r$
  \STATE $P_{\tilde c,\tilde r}\leftarrow 1$
  \STATE ${n_{\mathrm{iter}}} \leftarrow \min \left\{ {\left\lceil {\frac{{\log \left( {\frac{1}{\epsilon}} \right)}}{{2\log \left( {\frac{{{\mu _{\max }} + {\mu _{\min }}}}{{{\mu _{\max }} - {\mu _{\min }}}}} \right)}}} \right\rceil ,\eta } \right\} $
  \ELSE
  \STATE $\mathrm{count} \leftarrow \mathrm{count}+1$
  \ENDIF
  \ENDWHILE
  \end{algorithmic}
\end{algorithm}





 
 
\begin{figure}
    \centering
    \includegraphics[width=0.4\textwidth]{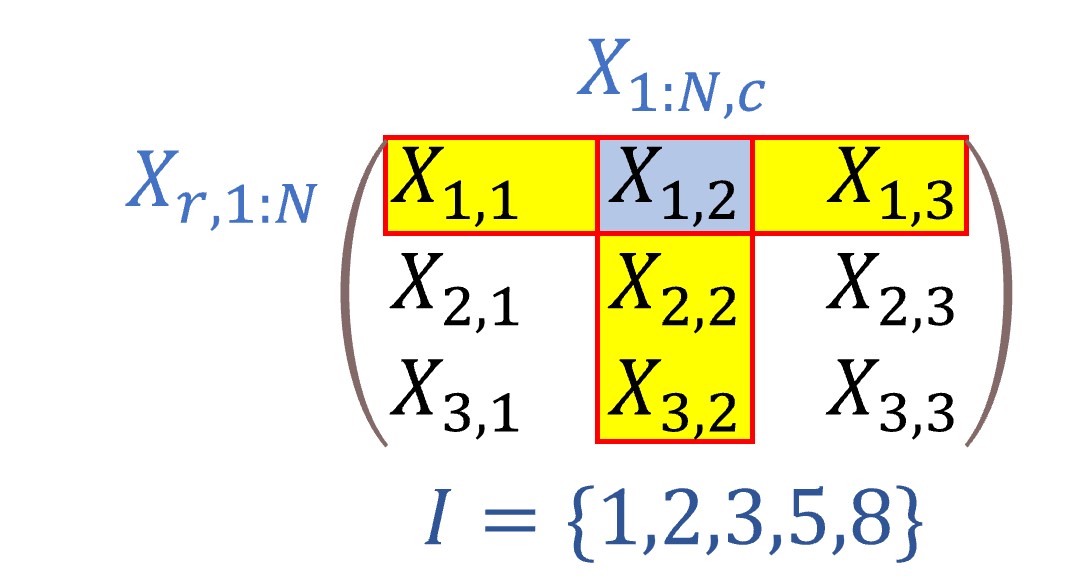}
    \caption{Deleting the $r$-th row and the $c$-th column from $X(k)$, where $X(k)_{r,c}$ is the largest element of $X(k)$.}
    \label{fig:reducex}
\end{figure}

\begin{figure}
    \centering
    \includegraphics[width=0.6\textwidth]{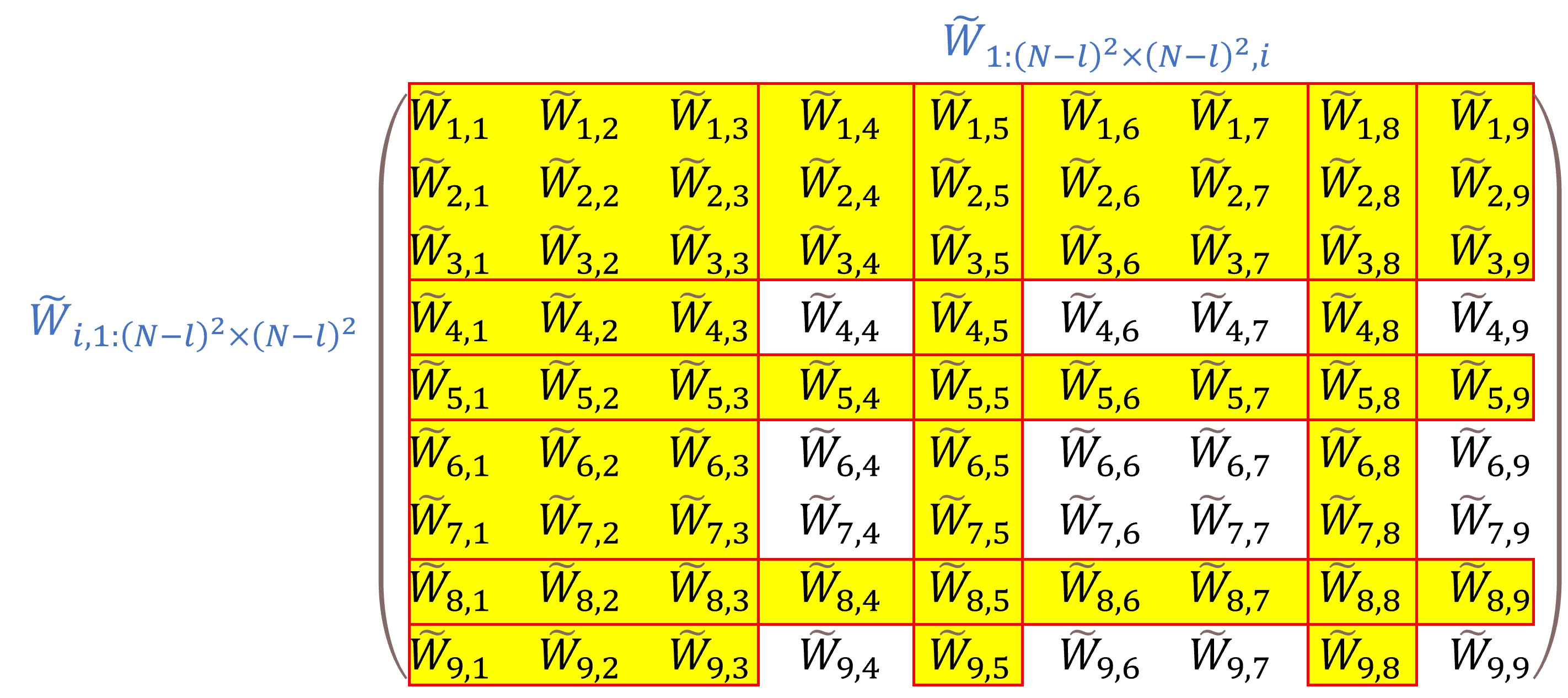}
    \caption{Deleting the $i$-th columns and rows from $\tilde W(l)$, where $i\in I$.}
    \label{fig:reducew}
\end{figure}

Line 2 is the steepest projection gradient descent operation, where $\nabla \tilde f(x(k))$ denotes the gradient of $\tilde f$ along $x(k)$, i.e. $\nabla \tilde f(x(k))=(\tilde W+\tilde W^T)x(k)+c$. As for the choice of time-varying gradient step size $\alpha_k$, we use the steepest descent algorithm in this paper, for it shows great convergence results for QP. In this algorithm the step size $\alpha_k$ is determined by means of optimal line search as $\alpha_k = \mathop {\arg \min }\limits_\alpha  \tilde f(x(k)) - \alpha \nabla \tilde f$, and for our quadratic cost function $\alpha_k$ has an explicit form solution $\alpha_k=\frac{{\nabla {\tilde f^T}({x_k})\nabla \tilde f({x_k})}}{{\nabla {\tilde f^T}({x_k})\tilde W\nabla \tilde f({x_k})}}$. Line 3 indicates time-triggering and we find the maximum elements of variable $X(k+1)$ in Line 6, indicated by column index $c$ and row index $r$. In Line 7, the corresponding rows and columns and deleted from the decision variable $X(k+1)$ (see Figure. \ref{fig:reducex}). By using $I$ to denote the index set of elements $X(k)_{r,1:N}$, $X(k)_{1:N,c}$ in $x(k)$, then $\tilde W(l)_{1:(N-l)^2\times (N-l)^2,i}$ and $\tilde W(l)_{i,1:(N-l)^2\times (N-l)^2}$ with $i\in I$ are deleted from the cost matrix $\tilde W(l)$ in Line 8 (see Figure \ref{fig:reducew}). The new cost matrix $\tilde W(l+1)$ has dimension of $(N-l-1)^2\times (N-l-1)^2$. 
After that, Line 9 aims at finding the original row and column index $\tilde c,\tilde r$, which are different from $c,r$. The reason is that $c,r$ belongs to the index set of the reduced matrix $X(k)$, with lower dimension compared with $X(0)$. Then, Line 10 sets element $P_{\tilde c,\tilde r}=1$ for the resulting permutation matrix $P$. Finally, Line 11 calculates the time coefficient for the next dimension reduction process depending on the predefined tolerance $\epsilon$, and a predefined positive integer $\eta$. Here, the formulation is derived from the linear convergence speed for steepest descent algorithms over quadratic cost \cite{meza2010steepest}. The exponential decay term is ${\left( {\frac{{{\mu _{\max}} - {\mu _{\min}}}}{{{\mu _{\max}} + {\mu _{\min}}}}} \right)^2}$, where $\mu_{\max}$ and $\mu_{\min}$ denote the maximum and minimum eigenvalues of the cost matrix $\tilde W(l)$, respectively. $\eta$ is used as a truncation to avoid $n_{\mathrm{iter}}$ becoming unacceptable large.

\subsection{Projection onto $\mathbb{DS}_N$}
 In Algorithm \ref{al1}, the projection function $Proj_{\mathbb{DS}_N}$ is defined as:
\begin{equation}\label{eq5}
\begin{array}{l}
Proj_{\mathbb{DS}_N}(X)=\mathop {\arg \min }\limits_{Y\in \mathbb{DS}_N} \frac{1}{2}||X - Y|{|_F^2},\\
\end{array}
\end{equation}
which yields the closest matrix in the set of doubly stochastic matrices based on Euclidean distance. Here we note that we omit an inverse vectorization operation in line 2 of Algorithm \ref{al1}. The authors of \cite{zass2006doubly} proposed an iterative algorithm for performing the projection onto doubly stochastic spaces, but only limited that to symmetric matrices. In most cases, the cost matrix can be symmetrized with no influence on the optimum. However, some special structure properties like sparsity may be changed. So here we give a more general alternating directional projection algorithm for random input matrix, where the doubly stochastic matrix $Y$ is not necessarily symmetric. This is summarized in Algorithm \ref{al2}.

\begin{algorithm}[h]
  \caption{Alternating Directional Projection Algorithm}
   \hspace*{\algorithmicindent} \textbf{Input:} matrix $X$\\
 \hspace*{\algorithmicindent} \textbf{Output:} doubly stochastic matrix Y\\
 \vspace{-3ex}
  \begin{algorithmic}[1]\label{al2}
  \WHILE{$Y \notin \mathbb{DS}_{N}$}
  \STATE $\lambda \leftarrow (AA^T)^{-1}(1-AX)$
  \STATE $X\leftarrow X+A^T\lambda/2 $
  \STATE $X\leftarrow th_{\ge 0}(X)$
  \ENDWHILE
  \STATE $Y \leftarrow X$
  \end{algorithmic}
\end{algorithm}

\begin{theorem}\label{th3}
The input matrix $X$ converges to the closest doubly stochastic approximation $P$ with Algorithm \ref{al2}.
\end{theorem}

\begin{proof}
We split the projection problem \eqref{eq5} into two sub-problems, where one is with an inequality constraint, and the other one is with an equality constraint. The advantage is that each of the sub-problems has an analytical solution. Therefore, the solution of \eqref{eq5} is at the intersection of the two sub-problems.

Consider the sub-problem with an affine equality constraint:
\begin{equation}\tag{$\text{P}_1$}
\begin{array}{l}
 \mathop {\arg \min }\limits_Y~~~~||X - Y||_F^2\\
\mathrm{subject~to}\;~AY=\bf{1},
\end{array}
\end{equation}
where $A = \left[ \begin{array}{l}
{I_N} \otimes {{\bf{1}}}^T\\
{{\bf{1}}}^T \otimes {I_N}
\end{array} \right]$. Note that the equality constraint for $\mathrm{vec}^{-1}(Y)\in \mathbb{DS}_N$ is:
\begin{equation*}
    \mathrm{vec}^{-1}(Y){\bf{1}} = {\bf{1}}, {\bf{1}}^T\mathrm{vec}^{-1}(Y) = {\bf{1}}^T.
\end{equation*}

Reformulating the above two equality constraints leads to $AY=\bf{1}$. 
The corresponding Lagrangian is
\begin{equation}\label{eq8}
    L(Y,\lambda ) = \text{Tr}({Y^T}Y - 2{X^T}Y) - 2 \lambda^T (AY - {\bf{1}}),
\end{equation}
where $\lambda$ is not constrained. The first order condition over the primal variable $Y$ of $L(Y,\lambda)$ results in:
\begin{equation}\label{eq9}
    Y=X+A^T\lambda.
\end{equation}
Multiplying \eqref{eq9} on both sides by $A$:
\begin{equation}\label{eq10}
    \begin{array}{l}
AY = AX + A{A^T}\lambda,
\end{array}
\end{equation}
implies that $\lambda=(AA^T)^{-1}({\bf 1}-AX)$. Combined with \eqref{eq9}, we get the explicit solution as:
\begin{equation}\label{eq11}
    X + {A^T}{(A{A^T})^{ - 1}}({\bf{1}} - AX).
\end{equation}

The second subproblem with affine inequality constraint is:
\begin{equation}\tag{$\text{P}_2$}
\begin{array}{l}
\mathop {\arg \min }\limits_Y ~~~~||X - Y||_F^2\\
\mathrm{subject~to}\;~Y \ge 0.
\end{array}
\end{equation}
Its solution is
\begin{equation}\label{eq12}
    \mathrm{th}_{\ge 0}(X).
\end{equation}

Then, with the iterative projection onto the two sets, Algorithm \ref{al2} leads to the convergence of $X$ to the projection onto the intersection, $\mathbb{DS}_N$.
\end{proof}

\subsection{Complexity Analysis}
This subsection presents the results of the complexity analysis of our algorithm, including two parts: the number of iterations required for convergence onto a permutation matrix and floating-point operations required. 
\begin{theorem}\label{th2}
The number of iterations required for the resulting matrix $P\in \Pi_N$ of Algorithm \ref{al2} is $\mathcal{O}(N)$, and the number of floating point operations is $\mathcal{O}(N^3)$.
\end{theorem}
\begin{proof}
One column and one row of $X$ are reduced every\\
$Nn_{\mathrm{iter}}\le N\left\lceil {\frac{{\log \left( {\frac{1}{\epsilon}} \right)}}{{2\log \left( {\frac{{{\mu _{\max}} + {\mu _{\min}}}}{{{\mu _{\max}} - {\mu _{\min}}}}} \right)}}} \right\rceil$
iterations, while the element of $P$ with the same corresponding column index $\tilde c$ and row index $\tilde r$ is set to $1$. The two indices, $\tilde c$ and $\tilde r$, are then deleted from the column index set $1,\ldots,N$ and row index set $1,\ldots,N$. This guarantees no repetition of a non-zero entry over every column and row of matrix $P$ and the existence of such entry. Besides, although $n_{\mathrm{iter}}$ changes across iterations, it is irrelevant to the dimension $N$. Thus, with $\mathcal{O}(N)$ iterations, the resulting matrix $P\in \Pi_N$.

The floating-point operations mainly lie in line 2 of Algorithm \ref{al1}. The total complexity of lines 3-13 can be omitted as it scales as $\mathcal{O}(N^2)$. We first analyze the complexity of gradient descent. The floating point operations required for $N^2$ dimensional $x(k)-\alpha_k\nabla f(x(k))$ is $\mathcal{O}(N^2)$, this operation repeats for $Nn_{\mathrm{iter}}$ times, therefore the total cost is $n_{\mathrm{iter}}\sum\nolimits_{n = 1}^N {{n^2}}=\mathcal{O}(n_{\mathrm{iter}}N^3)$. The floating-point operations required in Algorithm \ref{al2} mainly come from matrix multiplication, since the inverse calculation over $AA^T$ only needs to be done once. Thus, with $\mathcal{O}(N^2)$ time floating-point operations in a loop, the total amount of computation is $\mathcal{O}(N^3)$. Combining the complexity of gradient descent and projection, the resulting time complexity is $\mathcal{O}(N^3)$, as $n_{\mathrm{iter}}$ is a constant.
\end{proof}
\begin{remark}
The choice of $n_{\mathrm{iter}}$ depends on the accuracy and computational speed requirement. With larger $n_{\mathrm{iter}}$, the accuracy will be higher because of more gradient descent steps, whereas the computational speed will be lower.
\end{remark}

Compared to existing results on the convex-to-concave method \cite{fogel2013convex}, \cite{dym2017ds++}, \cite{bernard2018ds}, \cite{jiang2016l_p}, our method does not require to incrementally tune the penalty term to guarantee the solution to be a permutation matrix. In their formulations, the penalties often depend on the maximum and minimum eigenvalues of matrix $\tilde W$, which brings extra difficulties on computation for a large scale system. Besides, in our dimension reduction algorithm the accumulated time consumption of the multiplication is lower compared with the convex-to-concave method because of dimension reduction.

\subsection{A Note on Convexity}
It is proved that the relaxed QAP \eqref{RQAP} is convex if and only if $\tilde W$ is positive semi-definite, and strictly convex if and only if $\tilde W$ is positive definite. The following two problems are of interest:
\begin{problem}\label{pro1}
Can Algorithm 1 preserve the convexity of QAP in each iteration if it is convex initially?
\end{problem}

\begin{problem}\label{pro2}
Can Algorithm 1 preserve the $\sigma$-convexity of QAP in each iteration if it is  $\sigma$-strongly convex initially?
\end{problem}
Note that Problem \ref{pro1} is a special case of Problem \ref{pro2}. We prove the statements of Problem \ref{pro2}, which then directly extends to Problem \ref{pro1}.

\begin{theorem}\label{the:convex}
$\tilde f_l(x)$ is $\sigma$-strongly convex for all $l\in[1,\ldots,N-1]$ if $\tilde f_0(x)$ is $\sigma$-strongly convex.
\end{theorem}
\begin{proof}
Let $\tilde W(l)$ denote the $l$-th cost matrix corresponding to $\tilde f_l(x)$, then we have $\tilde W(l)\in \mathbb{R}^{(N-l)^2\times (N-l)^2}$. We use $\textbf{r}_l=\{r_l^1,\ldots,r_l^{2(N-l)-1}\}$ to represent the index set of reduced columns and rows from $\tilde W(l)$, $\textbf{h}_l=\{1,\ldots,N-l\} \backslash \textbf{r}_l$ to denote the indices set of residue elements. Then, we select $x(l)\in \mathbb{R}^{(N-l)^2}$ of which the $m$-th element  $x(l)_m=0, \forall m \in \textbf{r}_l$. We assume that $\tilde W(l)$ is $\sigma$-strongly convex. Then the following holds
\begin{equation}\label{eq15}
\begin{split}
        &{x(l)}^T(\tilde W(l) - \frac{\sigma}{2}I_{N-l}){x(l)} \ge 0,\\
        &\Rightarrow \sum\limits_{i \in {h_l}} {\sum\limits_{j \in {h_l}} {{x(l)_i}{x(l)_j}{\tilde W(l)_{i,j}}}  - \sum\limits_{i \in {h_l}} {\frac{{\sigma {x(l)_i}^2}}{2}} }  \ge 0.
\end{split}
\end{equation}
Let $x(l+1)$ denote the dimension reduced vector of $x(l)$, then \eqref{eq15} implies:
\begin{equation}\label{eq16}
    {x(l+1)}^T(\tilde W(l+1) - \frac{\sigma}{2}I_{N-l-1}){x(l+1)} \ge 0.
\end{equation}
Therefore, since $x(l)$ is defined randomly over indicies $\textbf{h}_l$, \eqref{eq16} proves that $\tilde W(l+1)$ is $\sigma$-strongly convex. In addition, it is assumed that $\tilde W(0)$ is $\sigma$-strongly convex, it follows that $\tilde W(l)-\frac{\sigma}{2}I_{N-l}$ is positive semi-definite for all $l\in[1,\ldots,N]$, which is equivalent to $f_l(x)$ is $\sigma$-strongly convex for all $l\in[1,\ldots,N]$. 
\end{proof}

\section{Simulation}
\label{sec:6}

We tested our algorithm on a variety of QAP instances, all the experimental data comes from QAPlib\footnote{\hyperlink{https://coral.ise.lehigh.edu/data-sets/qaplib/qaplib-problem-instances-and-solutions/}{https://coral.ise.lehigh.edu/data-sets/qaplib/qaplib-problem-instances-and-solutions/}} which includes Bur, Chr, Els, Esc, Had, Kra, Lipa, Nug, Rou, Scr, Sko, Ste, Tai, Tho, and Wil. We compare our results against two other algorithms, the SDP relaxation (C-SDP) \cite{ferreira2018semidefinite}, and the doubly stochastic relaxation \cite{fogel2013convex} with convex-to-concave \cite{zaslavskiy2008path} algorithm (PATH). These two algorithms stand for two kinds of convexfication algorithms shown in Section \ref{sec:4}. The former one belongs to SDP convexification algorithm \eqref{eq14} and the later one belongs to the doubly stochastic convexification algorithm \eqref{CRQAP}. The comparison includes two parts, i.e. the final value of cost and computing time. The dimension of data set varies from $10$ to $256$. It should be noted that the class of SDP relaxation algorithms \cite{de2015new}, \cite{zhao1998semidefinite}, \cite{ferreira2018semidefinite} is not able to produce permutation matrices. To illustrate the comparison, we utilize the Hungarian algorithm for projection onto permutation matrices, as introduced in the last paper. All the experiments are performed on a PC with 32GB RAM, 3.8GHz Intel i7-10700KF CPU.
\subsection{Computing Time}

\begin{figure}
     \centering

     \begin{subfigure}[b]{0.48\textwidth}
         \centering
         \includegraphics[width=0.8\textwidth]{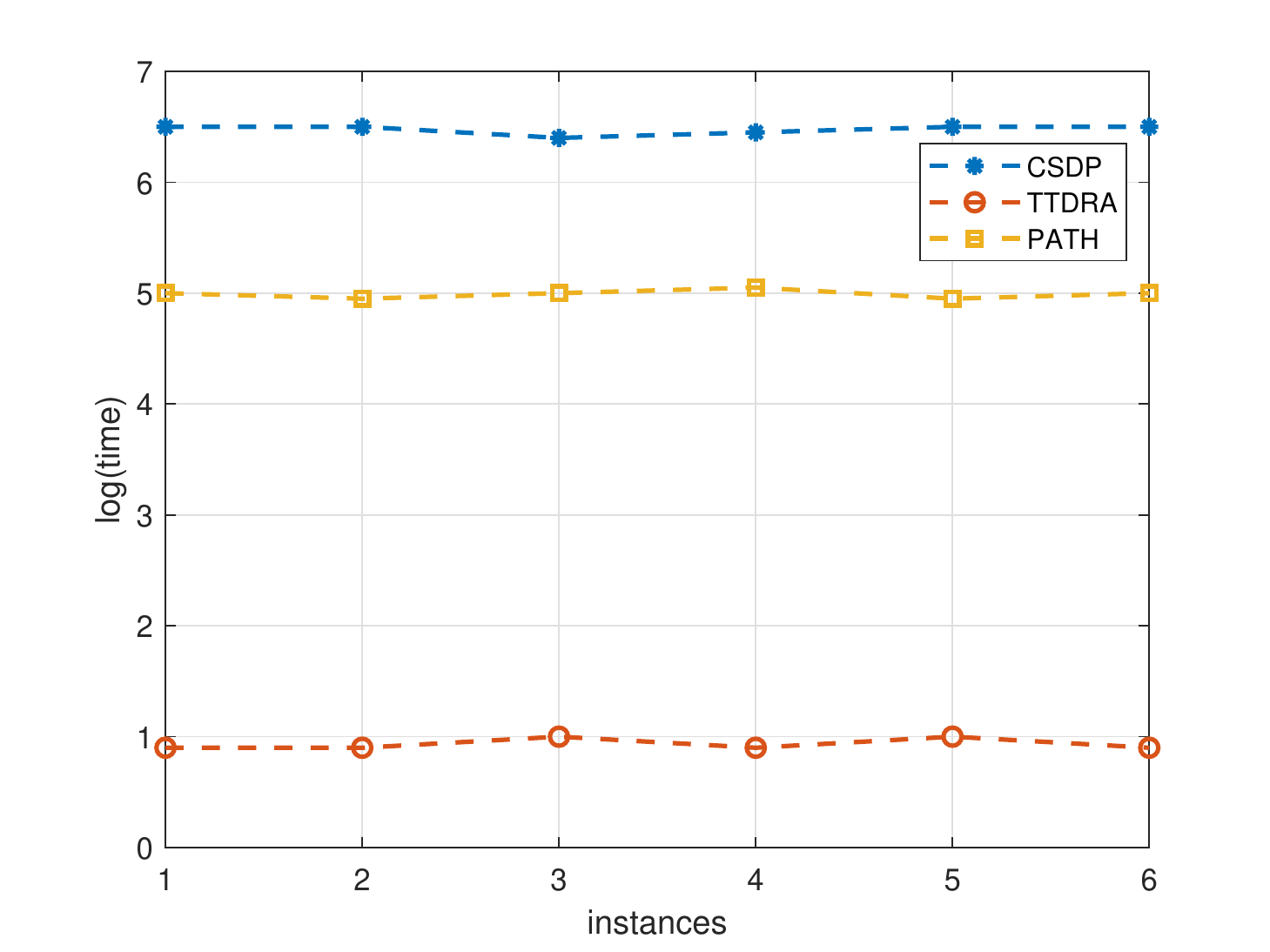}
         \caption{computation time for \textbf{bur} class instances}
         \label{fig:burtime}
     \end{subfigure}
     \hfill
     \begin{subfigure}[b]{0.48\textwidth}
         \centering
         \includegraphics[width=0.8\textwidth]{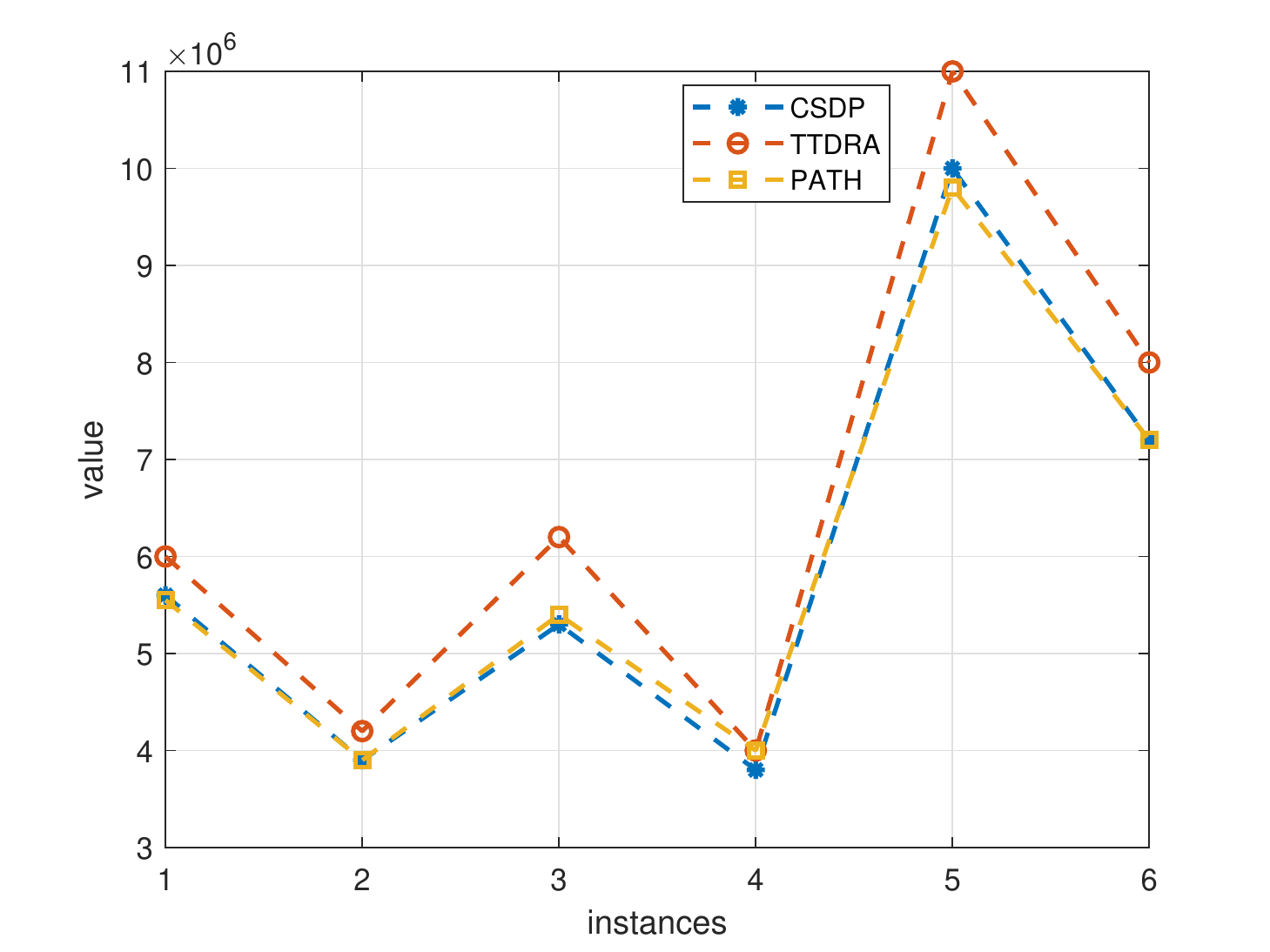}
         \caption{optimal value for \textbf{bur} class instances}
         \label{fig:burvalue}
     \end{subfigure}
     \hfill
     \begin{subfigure}[b]{0.48\textwidth}
         \centering
         \includegraphics[width=0.8\textwidth]{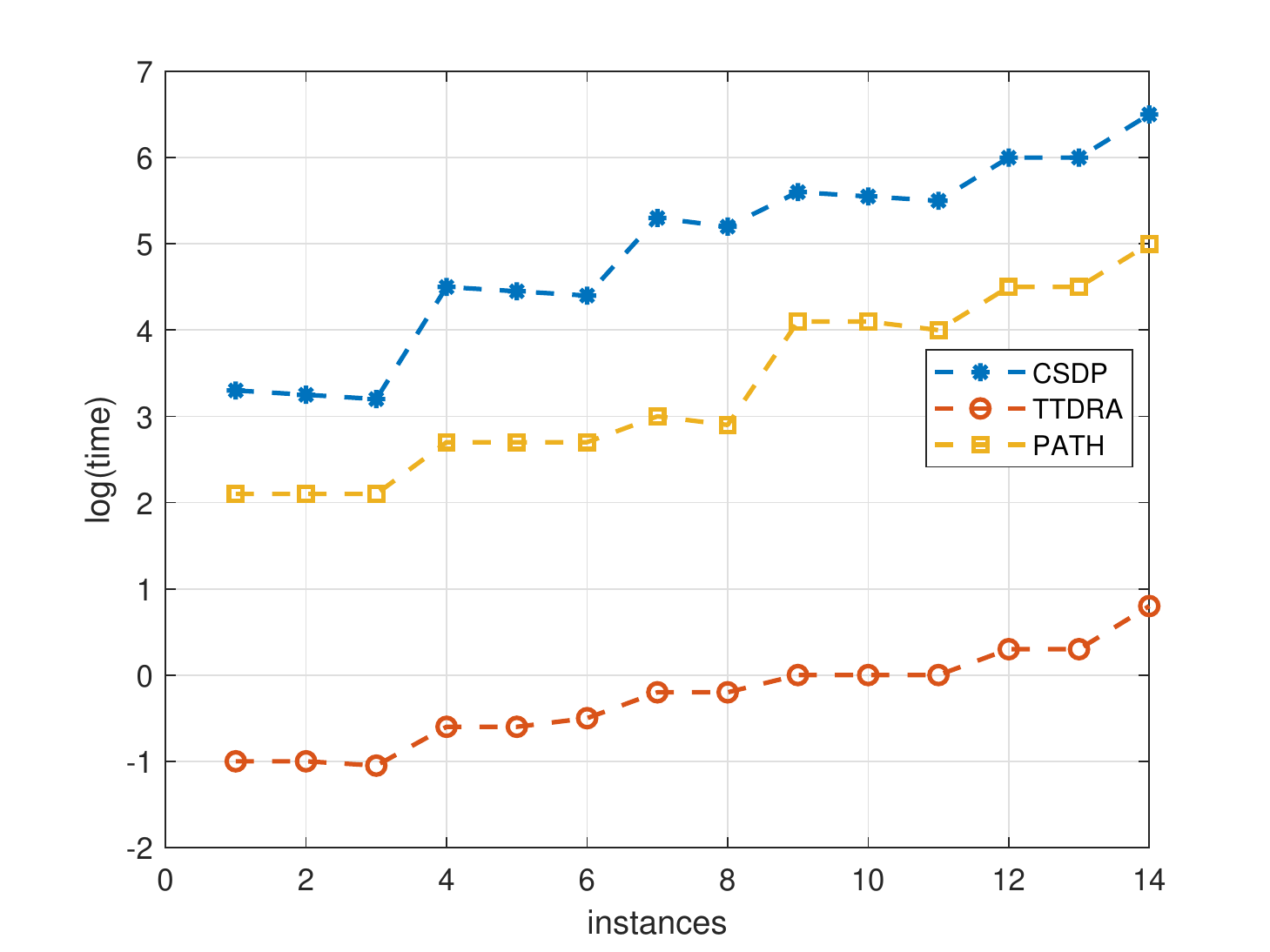}
         \caption{computation time for \textbf{chr} class instances}
         \label{fig:chrtime}
     \end{subfigure}
          \hfill
     \begin{subfigure}[b]{0.48\textwidth}
         \centering
         \includegraphics[width=0.8\textwidth]{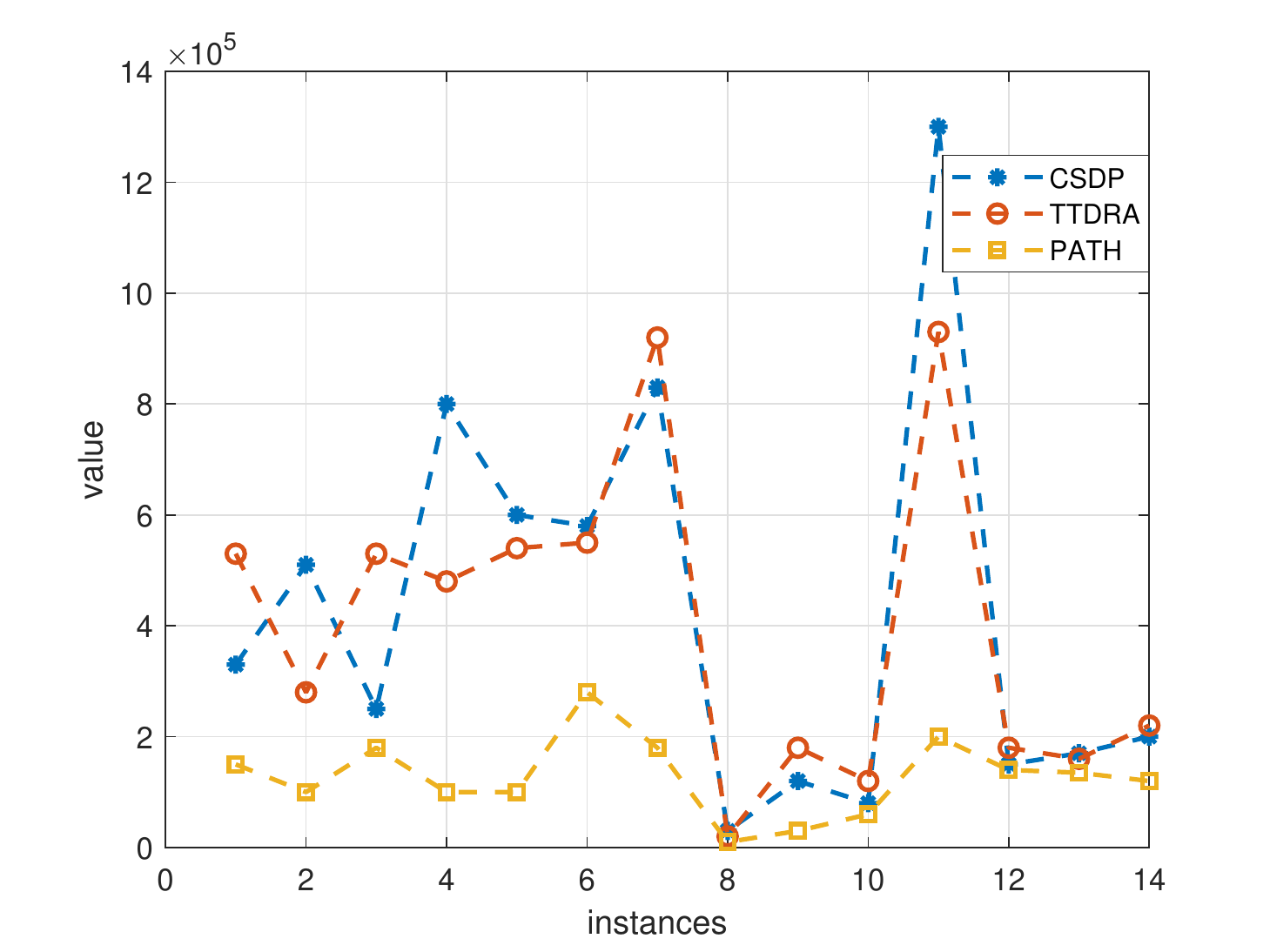}
         \caption{optimal value for \textbf{chr} class instances}
         \label{fig:chrvalue}
     \end{subfigure}
          \hfill
     \begin{subfigure}[b]{0.48\textwidth}
         \centering
         \includegraphics[width=0.8\textwidth]{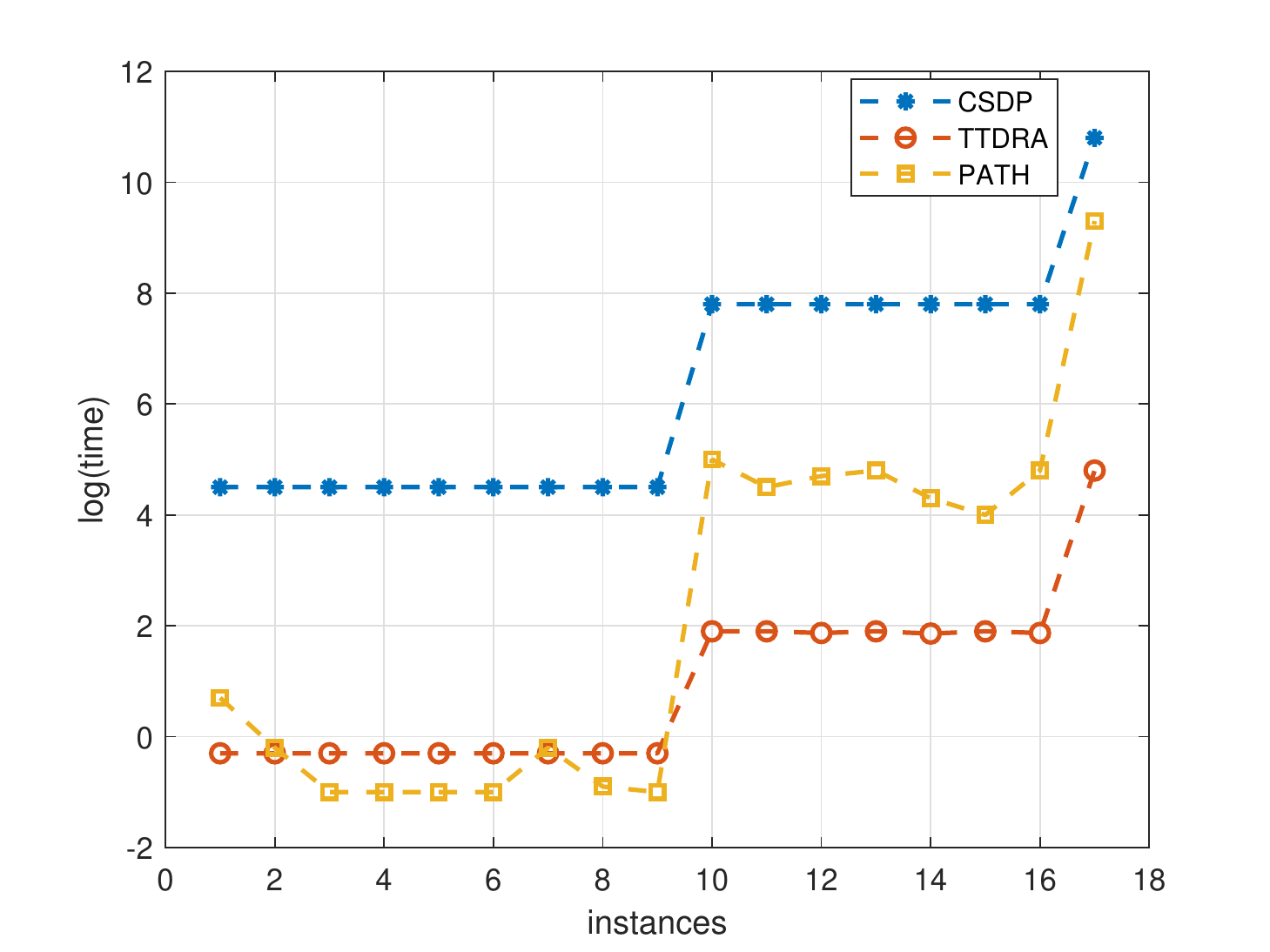}
         \caption{computation time for \textbf{esc} class instances}
         \label{fig:esctime}
     \end{subfigure}
          \hfill
     \begin{subfigure}[b]{0.48\textwidth}
         \centering
         \includegraphics[width=0.8\textwidth]{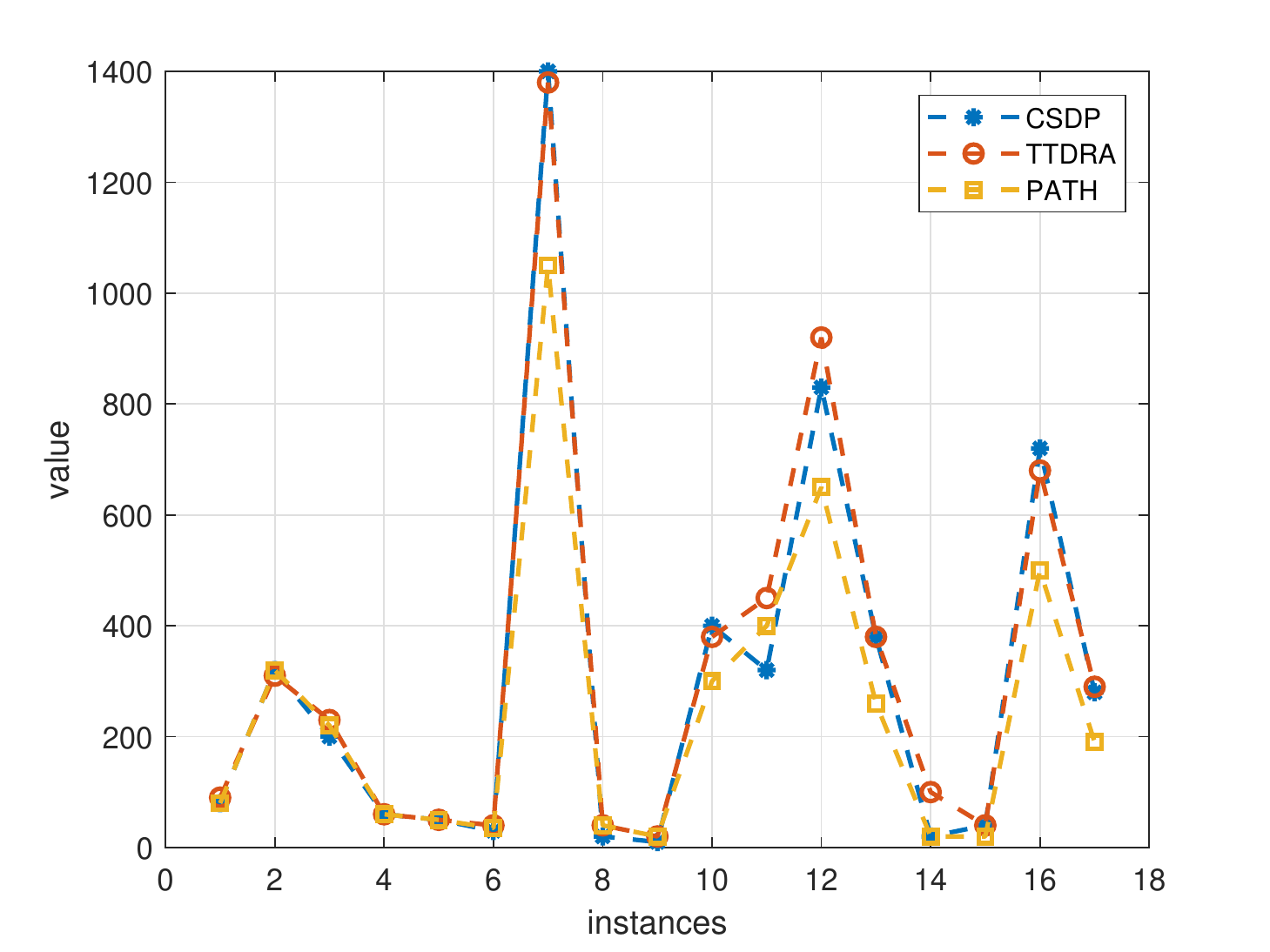}
         \caption{optimal value for \textbf{esc} class instances}
         \label{fig:escvalue}
     \end{subfigure}
          \hfill
     \begin{subfigure}[b]{0.48\textwidth}
         \centering
         \includegraphics[width=0.8\textwidth]{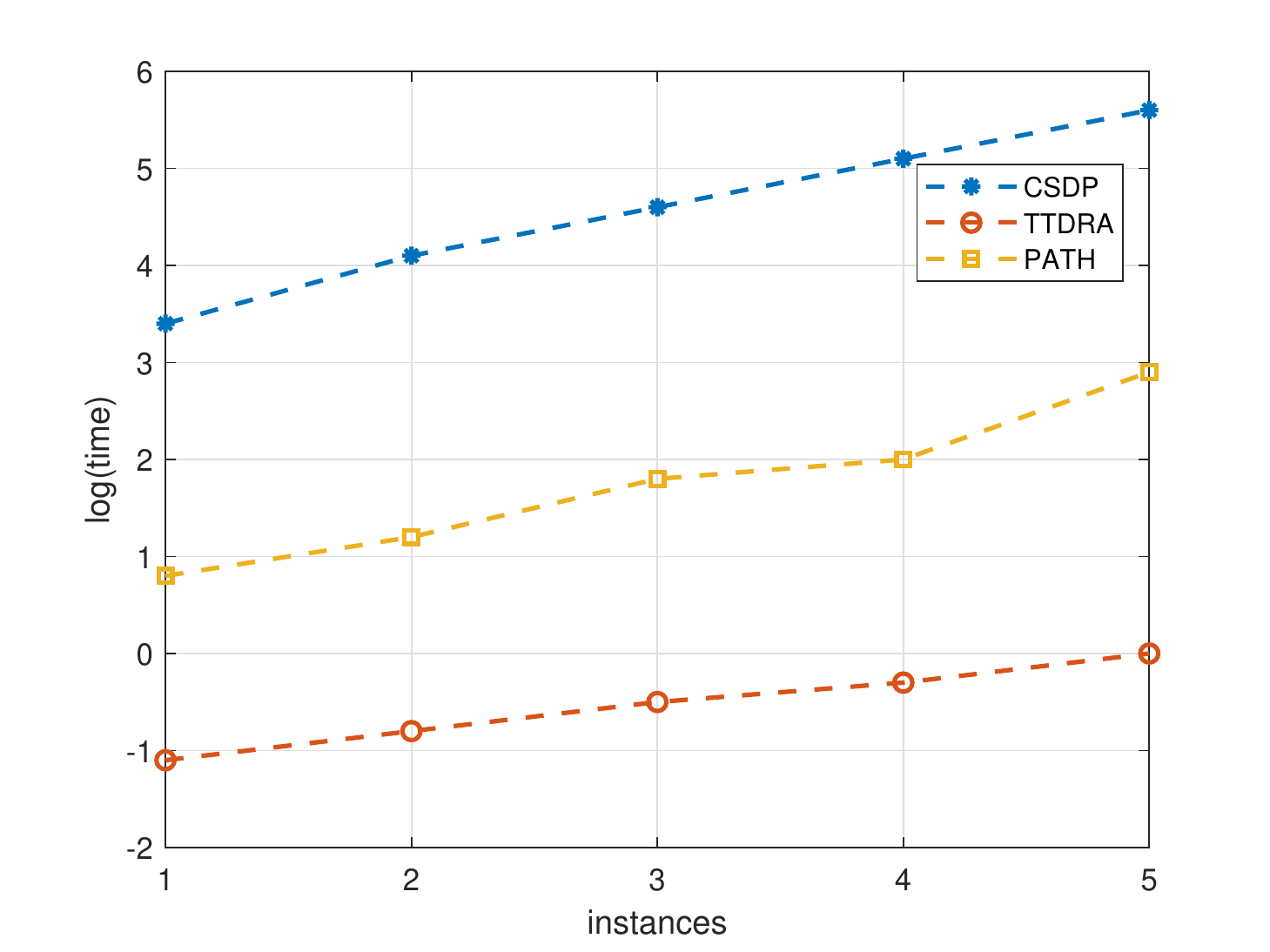}
         \caption{computation time for \textbf{had} class instances}
         \label{fig:hadtime}
     \end{subfigure}
              \hfill
     \begin{subfigure}[b]{0.48\textwidth}
         \centering
         \includegraphics[width=0.8\textwidth]{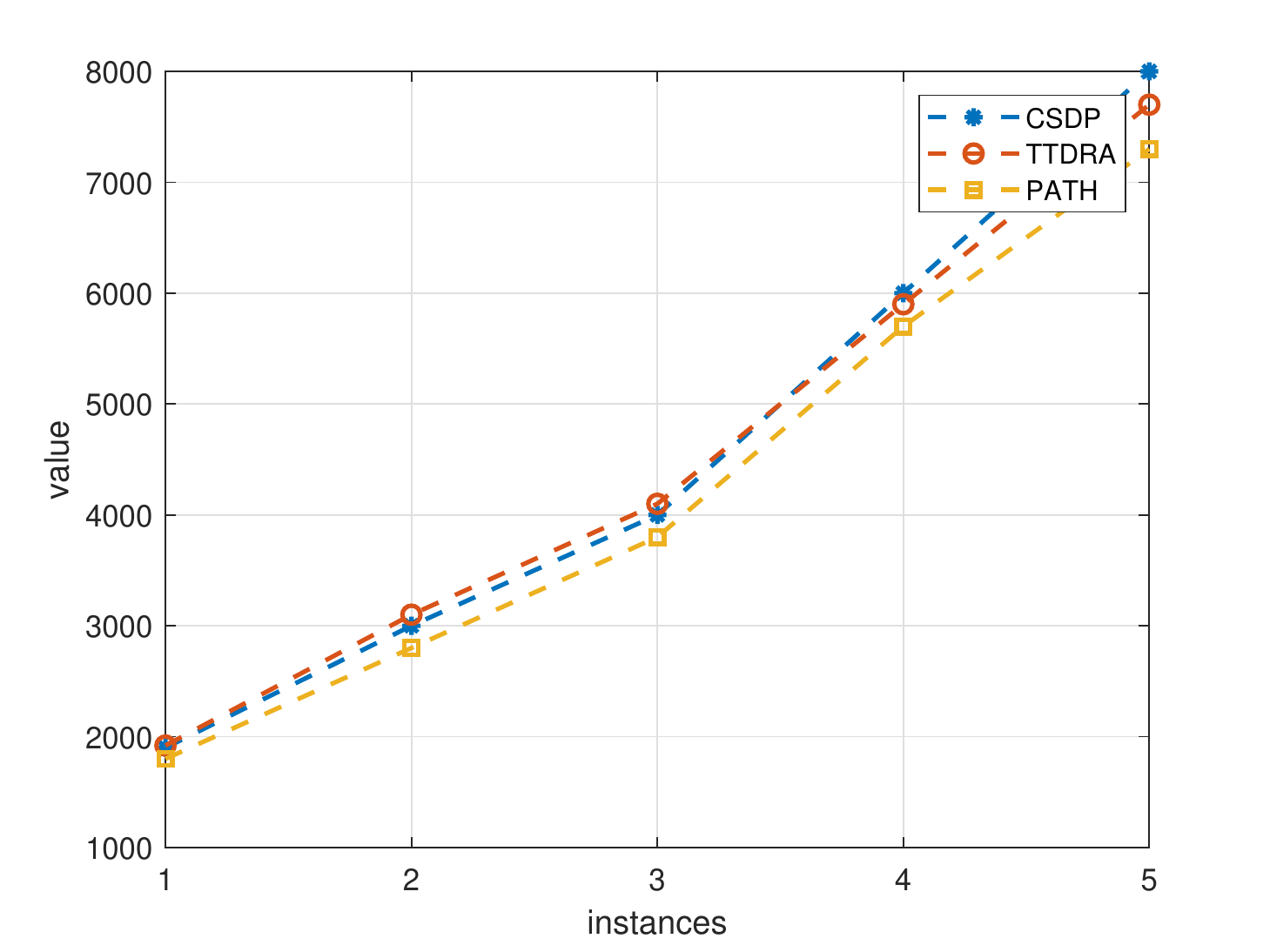}
         \caption{optimal value for \textbf{had} class instances}
         \label{fig:hadvalue}
     \end{subfigure}
         \caption{comparison between TTDRA, PATH, and C-SDP}
         \label{fig:timecomparison}
\end{figure}

We first show the optimality comparison among C-SDP, TTDRA, PATH. For C-SDP, we select $n=4$ graph nodes per variable, for PATH we chose $n=10$ iterations for convex-to-concave sampling, TTDRA we set the tolerance to be 0.5, and the convexity to be $10^6$. The first column of Fig. \ref{fig:timecomparison} shows the computation time used by TTDRA, PATH, and C-SDP. TTDRA is shown to be $10-10^4$ faster than PATH and $10^4-10^6$ faster than C-SDP; the C-SDP is the slowest algorithm since it lifts the dimension of the decision vector and requires solving an SDP, which is known to be computationally expensive. For some examples belonging to esc class, the computing speed of TTDRA is slower than PATH (Fig. \ref{fig:esctime}). This phenomenon is raised because of redundant iterations used for TTDRA, i.e. the minimum has been reached before the time condition is triggered.
    
\subsection{Optimality}
We then show the optimality comparison among C-SDP, PATH, and TTDRA. The second column of Fig. \ref{fig:timecomparison} shows the optimal values obtained by TTDRA, PATH, and C-SDP. We can see that the value obtained by TTDRA is competitive compared with other methods. On the other hand, the optimal value of TTDRA and C-SDP are close in most instances, and PATH acquired a better solution for more cases. The reason behind this results is that, PATH use the results from the last iteration.

\begin{figure}         
     \centering

     \begin{subfigure}[b]{0.48\textwidth}
         \centering
         \includegraphics[width=0.8\textwidth]{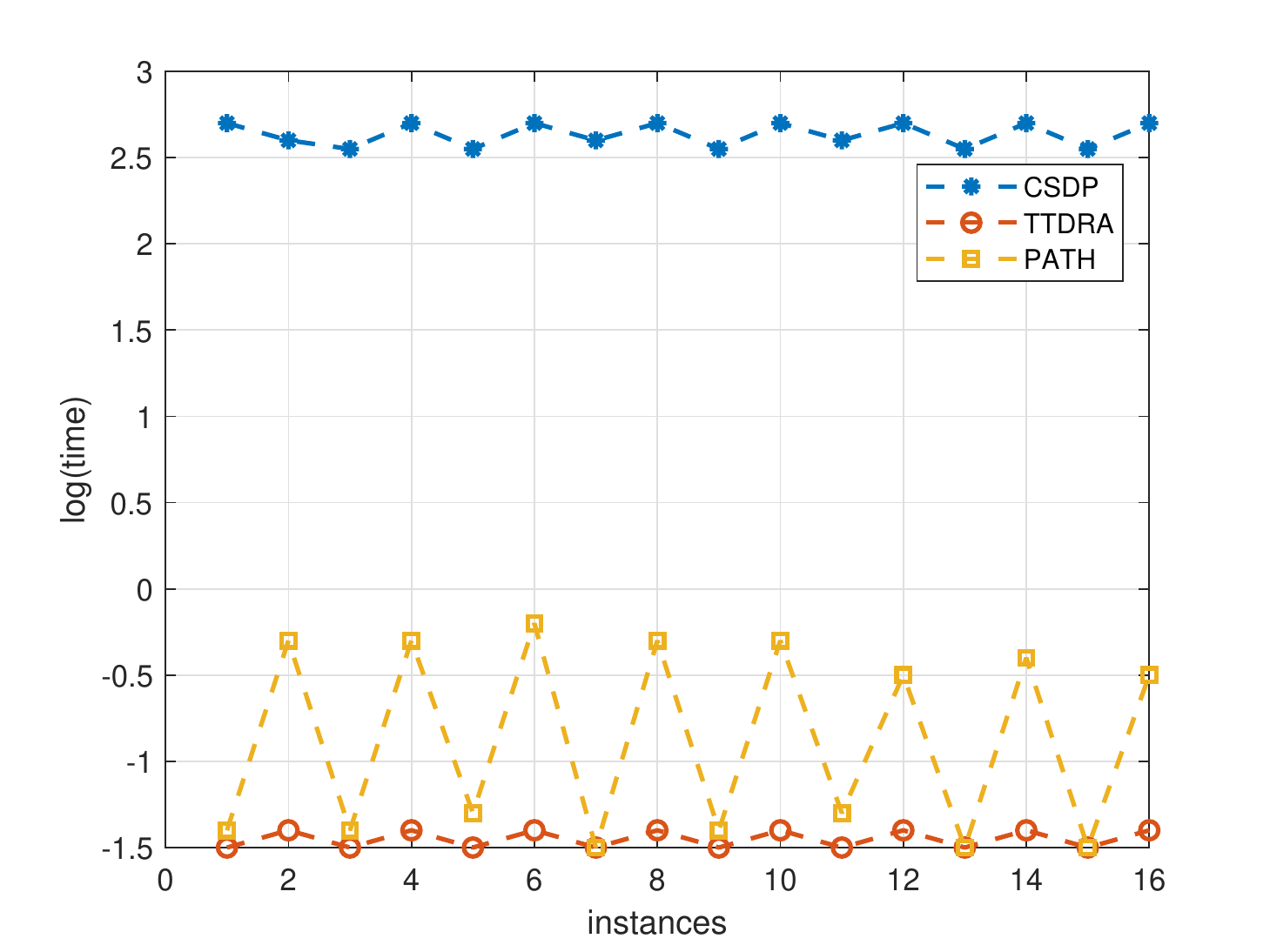}
         \caption{computational time for \textbf{lipa} class instances}
         \label{fig:lipatime}
     \end{subfigure}
     \hfill
     \begin{subfigure}[b]{0.48\textwidth}
         \centering
         \includegraphics[width=0.8\textwidth]{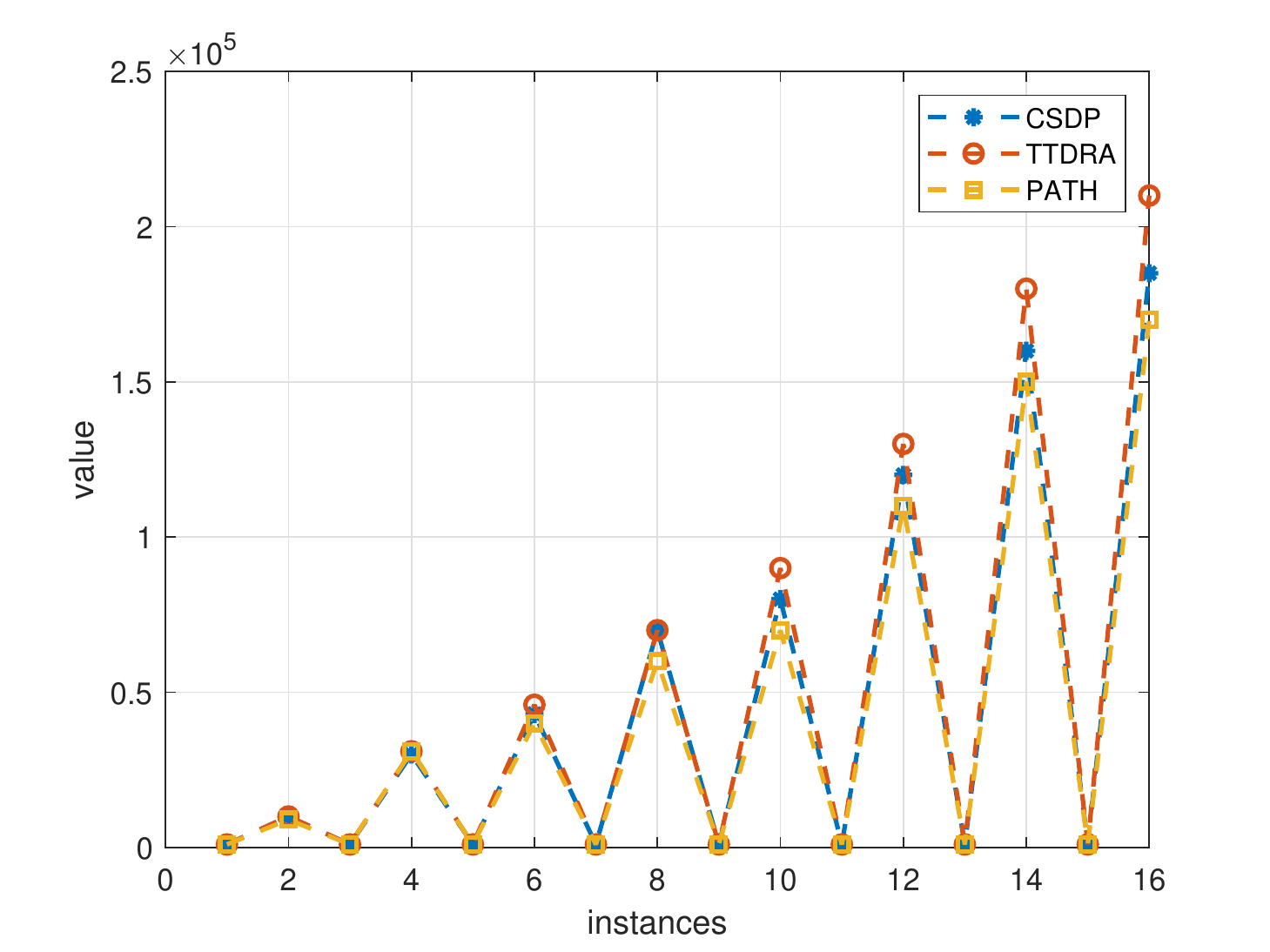}
         \caption{optimal value for \textbf{lipa} class instances}
         \label{fig:lipavalue}
     \end{subfigure}
     \hfill
     \begin{subfigure}[b]{0.48\textwidth}
         \centering
         \includegraphics[width=0.8\textwidth]{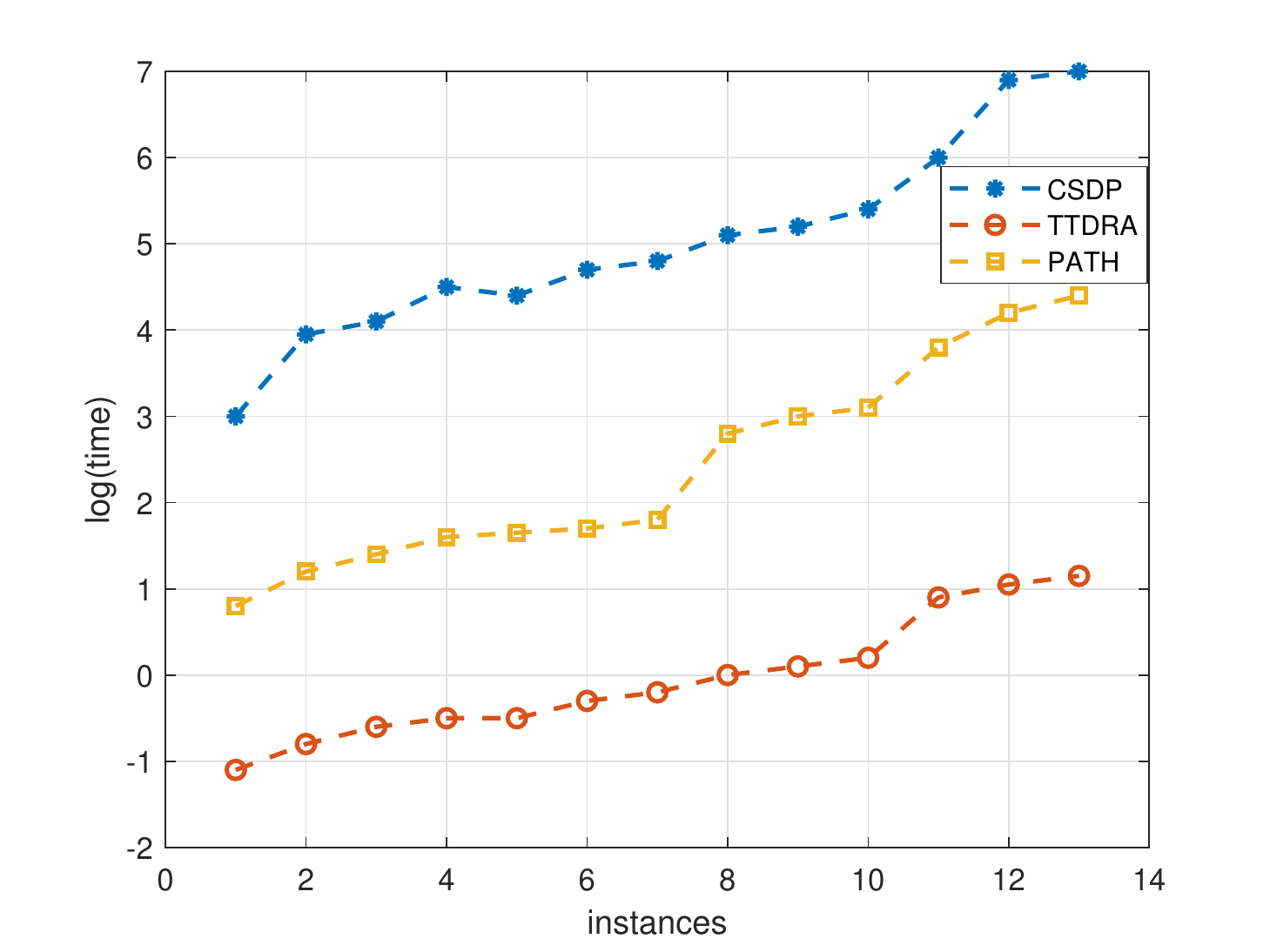}
         \caption{computational time for \textbf{nug} class instances}
         \label{fig:nugtime}
     \end{subfigure}
          \hfill
     \begin{subfigure}[b]{0.48\textwidth}
         \centering
         \includegraphics[width=0.8\textwidth]{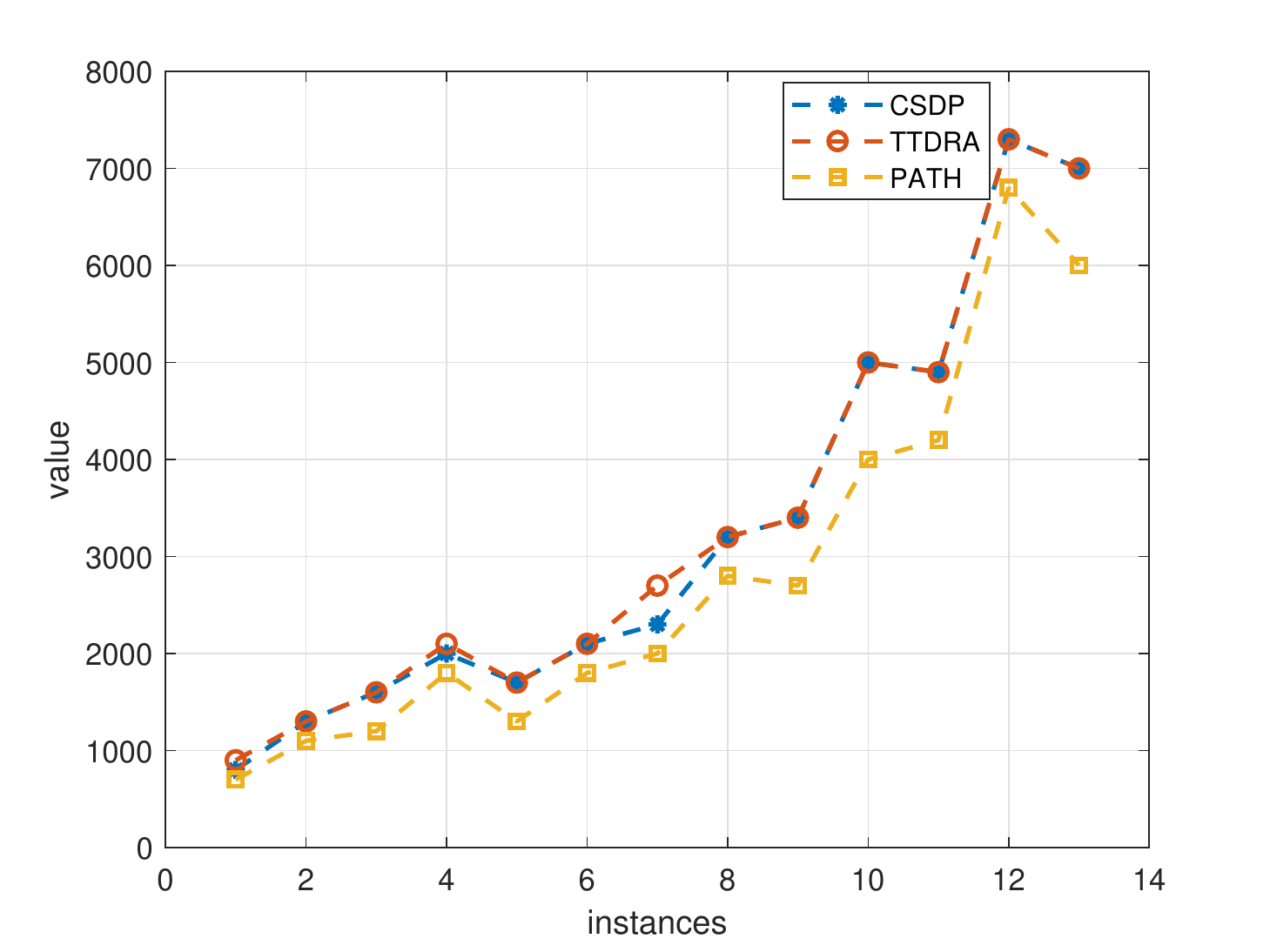}
         \caption{optimal value for \textbf{nug} class instances}
         \label{fig:nugvalue}
     \end{subfigure}
          \hfill
     \begin{subfigure}[b]{0.48\textwidth}
         \centering
         \includegraphics[width=0.8\textwidth]{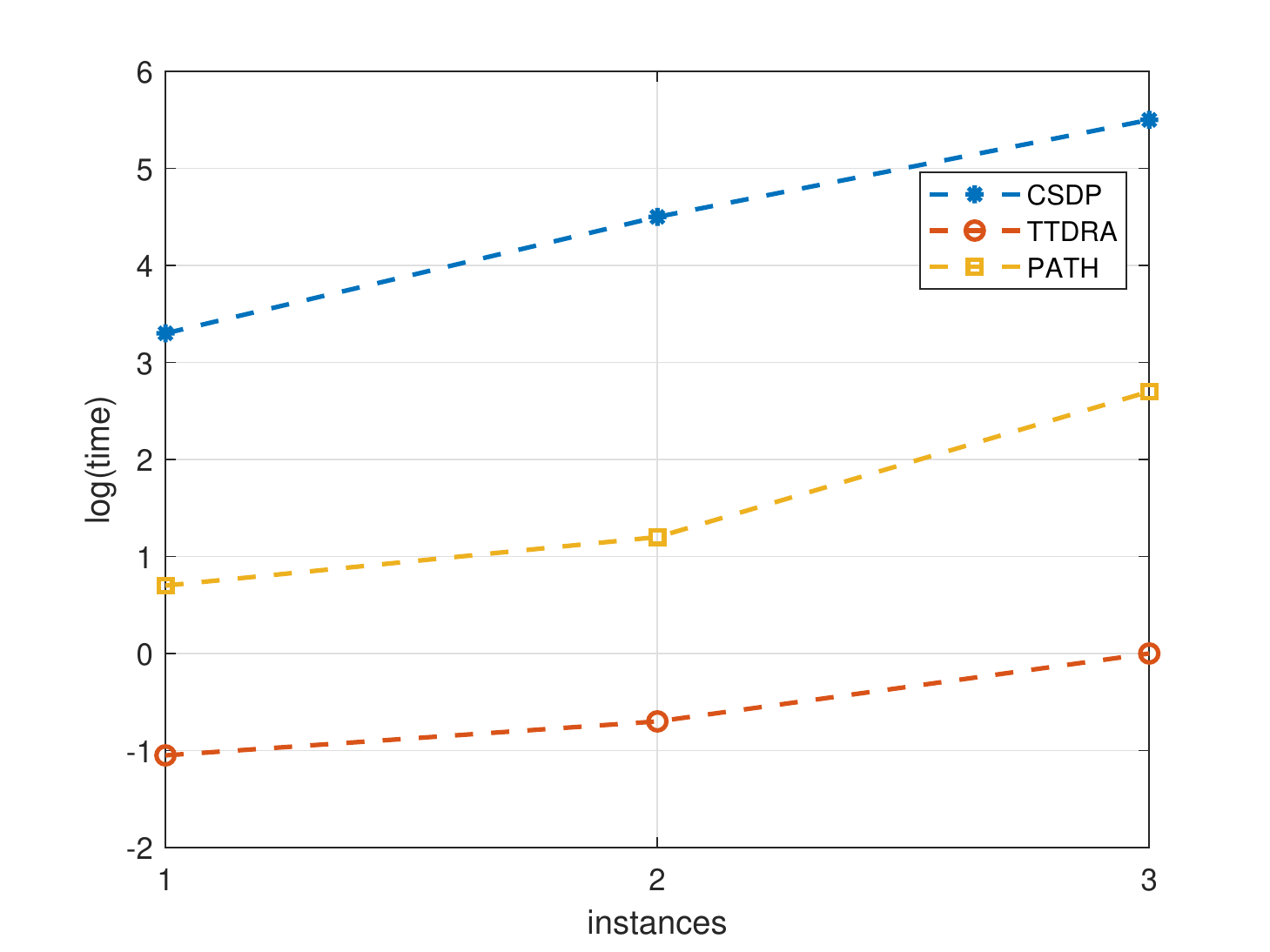}
         \caption{computational time for \textbf{rou} class instances}
         \label{fig:routime}
     \end{subfigure}
          \hfill
     \begin{subfigure}[b]{0.48\textwidth}
         \centering
         \includegraphics[width=0.8\textwidth]{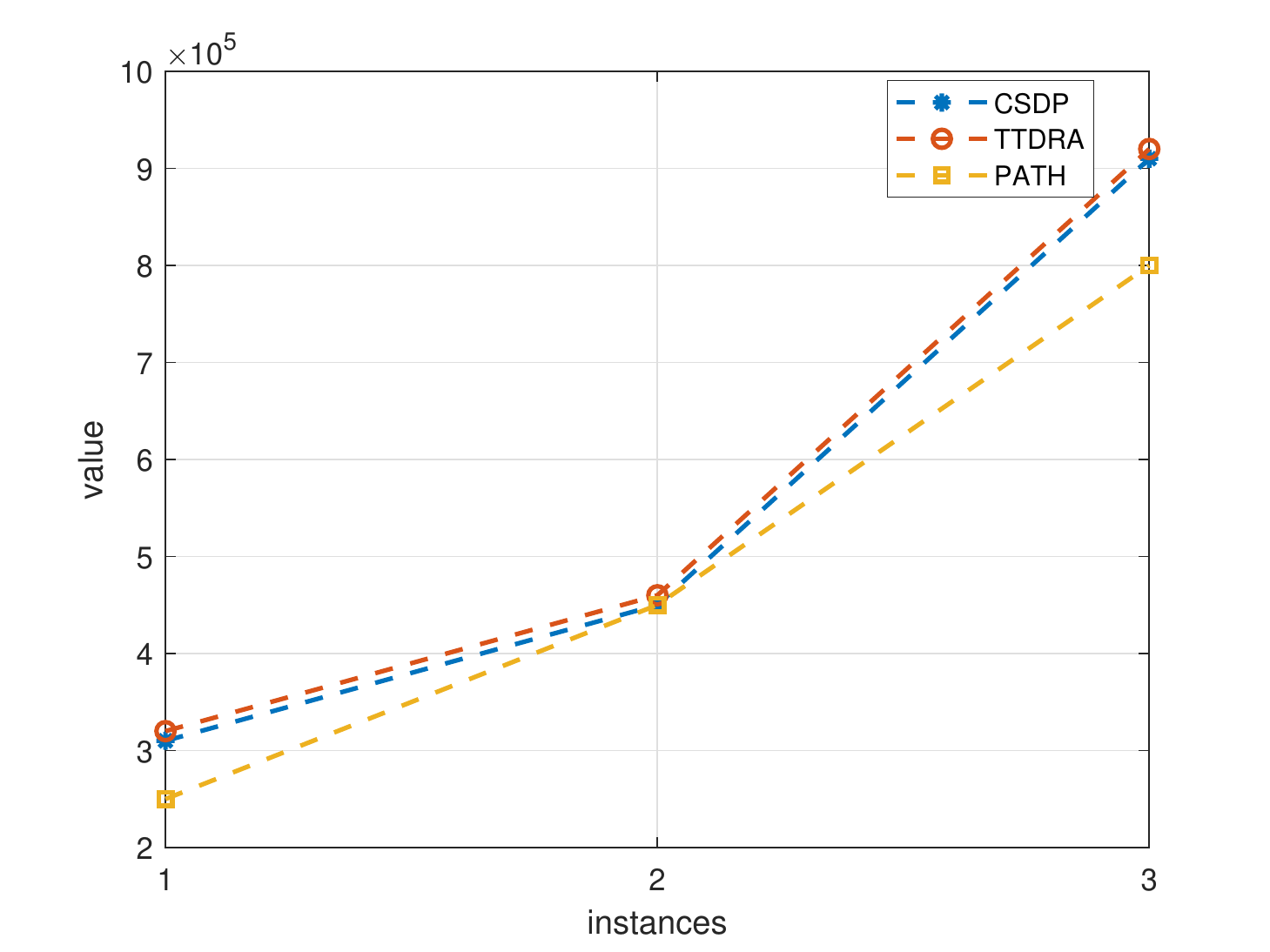}
         \caption{optimal value for \textbf{rou} class instances}
         \label{fig:rouvalue}
     \end{subfigure}
          \hfill
     \begin{subfigure}[b]{0.48\textwidth}
         \centering
         \includegraphics[width=0.8\textwidth]{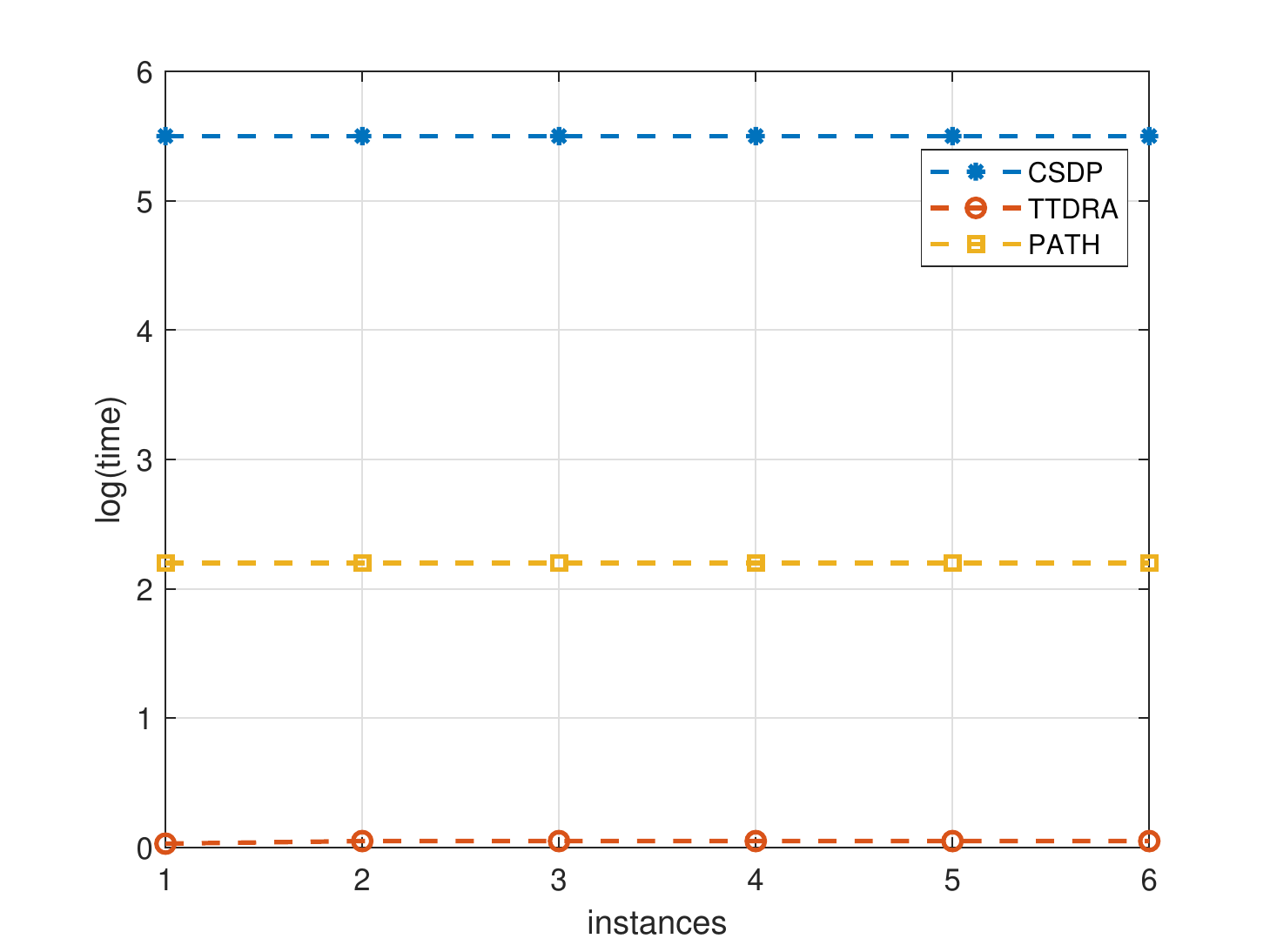}
         \caption{computational time for \textbf{sko} class instances}
         \label{fig:skotime}
     \end{subfigure}
              \hfill
     \begin{subfigure}[b]{0.48\textwidth}
         \centering
         \includegraphics[width=0.8\textwidth]{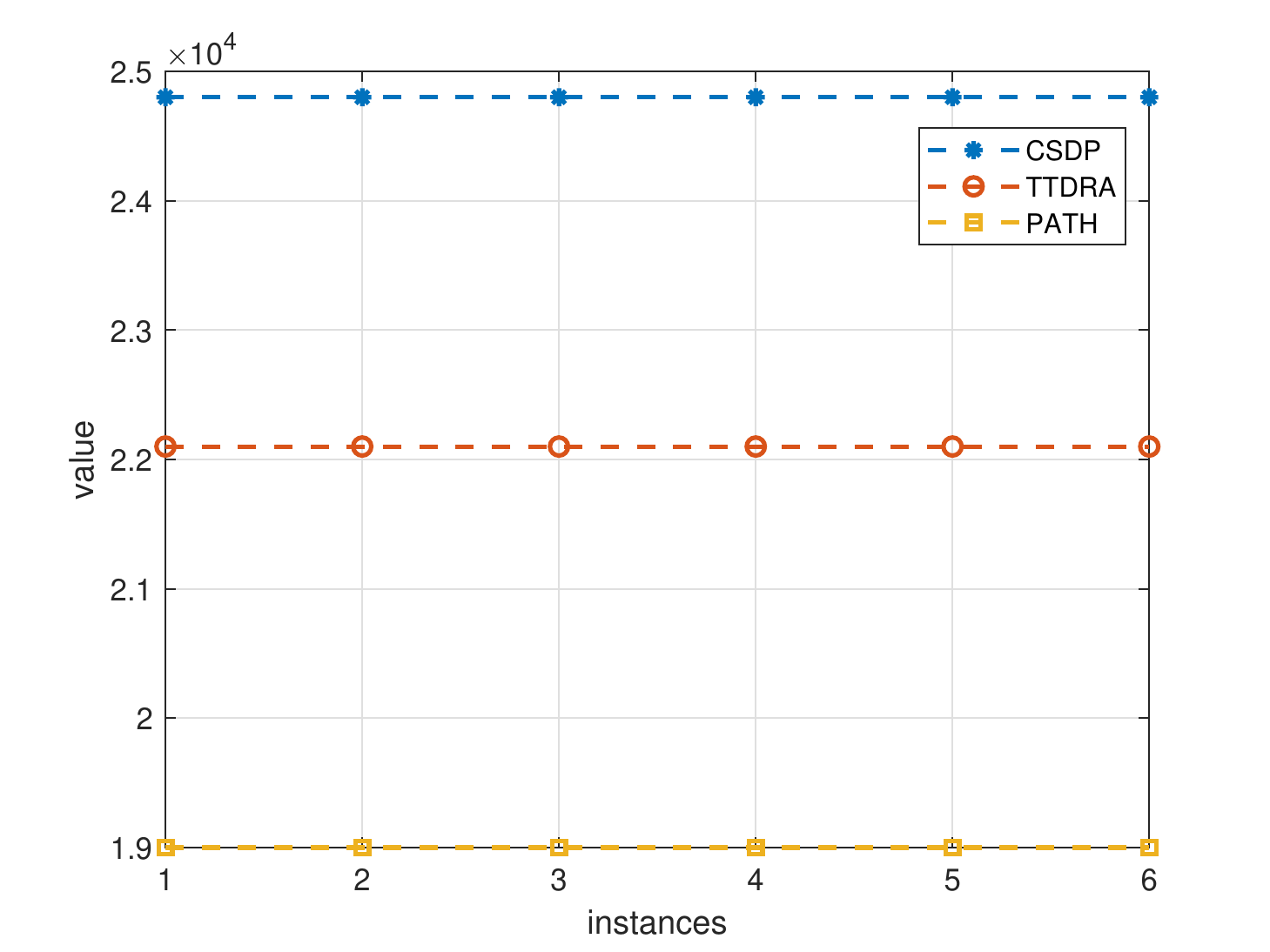}
         \caption{optimal value for \textbf{sko} class instances}
         \label{fig:skotime}
     \end{subfigure}
         \caption{comparison between TTDRA, PATH, and C-SDP}
\label{fig:valuecomparison}
\end{figure}

\section{Conclusion}
\label{sec:conclusion}
We have presented a time triggered dimension reduction algorithm for efficiently solving the task assignment problem. The nonconvex optimisation problem is convexified to be $\sigma$-strongly convex. The output of the algorithm is guaranteed to be a permutation matrix. We further showed that the convexity is preserved across the iterations. We also gave an upper bound of the computational complexity. The computational speed and optimality of our algorithm are verified on benchmark examples. In the future we aim at investigating stochastic variants of the proposed scheme, as well as parallelizable algorithms. 
\bibliographystyle{elsarticle-num} 
\bibliography{article.bib} 




\end{document}